\documentclass[aos,MSNbibl,nameyear,seceqn,dvips]{arximspdf}

\doi{10.1214/14-AOS1263} 
\volume{43}
\issue{1}
\pubyear{2015}
\firstpage{30}
\lastpage{56}
\docsubty{FLA}

\makeatletter

\newtheorem{thmm}{Theorem}[section]
\newtheorem{corollary}[thmm]{Corollary}
\newtheorem{lemma}[thmm]{Lemma}
\newproclaim{definition}{Definition}[section]
\newproclaim{remark}{Remark}[section]
\newproclaim{assumption}{Assumption}

\def\bfa{\mathbf a}
\def\bfc{\mathbf c}
\def\bfe{\mathbf e}
\def\bff{\mathbf f}
\def\bfx{\mathbf x}

\def\bfC{\mathbf C}
\def\bfD{\mathbf D}
\def\bfF{\mathbf F}
\def\bfI{\mathbf I}
\def\bfK{\mathbf K}
\def\bfM{\mathbf M}
\def\bfP{\mathbf P}
\def\bfS{\mathbf S}
\def\bfW{\mathbf W}
\def\bfX{\mathbf X}
\def\bfY{\mathbf Y}
\def\bfZ{\mathbf Z}

\def\bfomega{\bolds\omega}
\def\bftheta{\bolds\theta}
\def\bfbeta{\bolds\beta}
\makeatother

\begin{document}
\begin{frontmatter}

\title{Saturated locally optimal designs under differentiable
optimality criteria}
\runtitle{Saturated optimal designs under differentiable criteria}

\begin{aug}
\author[A]{\fnms{Linwei}~\snm{Hu}\thanksref{T1}\ead[label=e1]{pangaea@uga.edu}},
\author[B]{\fnms{Min}~\snm{Yang}\thanksref{T2}\ead[label=e2]{myang2@uic.edu}}
\and
\author[C]{\fnms{John}~\snm{Stufken}\corref{}\thanksref{T3}\ead[label=e3]{jstufken@asu.edu}}

\thankstext{T1}{Supported in part by NSF Grant DMS-10-07507.}
\thankstext{T2}{Supported in part by NSF Grants DMS-13-22797 and DMS-14-07518.}
\thankstext{T3}{Supported in part by NSF Grants DMS-10-07507 and DMS-14-06760.}
\runauthor{L. Hu, M. Yang and J. Stufken}

\affiliation{University of Georgia,
University of Illinois at Chicago
and Arizona State University}

\address[A]{L. Hu\\
Department of Statistics\\
University of Georgia\\
Athens, Georgia 30602\\
USA\\
\printead{e1}}

\address[B]{M. Yang\\
Department of Mathematics, Statistics\\
\quad and Computer Sciences\\
University of Illinois at Chicago\\
Chicago, Illinois 60607\\
USA\\
\printead{e2}}

\address[C]{J. Stufken\\
School of Mathematical\\
\quad and Statistical Sciences\\
Arizona State University\\
Tempe, Arizona 85287\\
USA\\
\printead{e3}}
\end{aug}
%

\received{\smonth{3} \syear{2014}}

%
\begin{abstract}
We develop general theory for finding locally optimal designs in a class
of single-covariate models under any differentiable optimality
criterion.
Yang and Stufken
[\textit{Ann. Statist.} \textbf{40} (2012) 1665--1681]
and
Dette and Schorning
[\textit{Ann. Statist.} \textbf{41} (2013) 1260--1267] gave
complete class results for optimal designs under such models. Based on
their results, saturated optimal designs exist; however, how to find
such designs has not been addressed. We develop tools to find
saturated optimal designs, and also prove their uniqueness under mild
conditions.
\end{abstract}

%
\end{frontmatter}
%
\section{Introduction}\label{sec:Introduction}
We consider the problem of finding locally optimal designs for a class
of single-covariate models under differentiable optimality criteria. In
order to avoid intricacies caused by the discreteness of the problem,
we will work with \emph{approximate} designs (see Section~\ref{sec:locally}). Because
the information matrix usually depends on the unknown parameters, we
consider \emph{locally} optimal designs by plugging in values for the
parameters in the information matrix. This gives good designs when
prior knowledge of the parameters is available, and it also provides a
benchmark for evaluating other designs. For the sake of simplicity, we
omit the word locally hereafter.

We provide general theoretical results that help to find saturated
optimal designs for many of the models for which previous results, such
as in \citet{ys12} and \citet{ds13}, have
established so-called \emph{complete class results}. While efficient
numerical algorithms, even without using the complete class results,
can be developed to approximate optimal designs, theory provides
unified results, both with respect to models and optimality criteria,
and offers insights that cannot be obtained from algorithms. In some
instances the theory enables us to find closed-form optimal designs;
moreover, it can be used to develop faster and better algorithms. For
example, because of the theory we can avoid having to discretize the
design space. We also use the theory to develop uniqueness results
under mild conditions, which cannot be obtained from an algorithm approach.

Our work is based on the complete class results given in a series of
papers, including most recently \citet{ys12} and \citet{ds13}. Based on their results, optimal designs can be found
in a small class of designs called the complete class, and in many
cases, this complete class only contains designs with at most $d$
design points, where $d$ is the number of parameters. However, theory
and tools to identify optimal designs for multiple optimality criteria
within the complete class have not been developed. So in Section~\ref{sec:locally}, we will present theorems to find optimal designs in
these classes and prove their uniqueness. Section~\ref{sec:app} applies
the theorems to a variety of different models including polynomial
regression models, nonlinear regression models and generalized linear
models. The computational benefits will be shown in Section~\ref{sec:computation}. Finally, Section~\ref{sec:Disscusion} gives a short
discussion about limitations of the approach. The technical proofs have
been relegated to the \hyperref[sec:Proofs]{Appendix}.

\section{Locally optimal design}
\label{sec:locally}
The models under consideration include polynomial regression models,
nonlinear regression models and generalized linear models, with a
univariate response $y$ and a single covariate $x$ which belongs to the
design space $[L,U]$ ($L$ or $U$ could be $-\infty$ or $\infty$,
resp., with $[L,U]$ being half open or open). The unknown
parameter is a $d\times1$ vector denoted as $\bftheta=(\theta_1,\ldots
,\theta_d)^T$. To be specific, for polynomial regression models and
nonlinear models $\bftheta$ is the unknown parameter in the mean
response $\eta(x,\bftheta)=\mathrm{E}(y)$. We assume the variance to be
constant unless otherwise specified, and take its value to be 1 since
it does not affect the optimal design. For generalized linear models,
$\bftheta$ is the unknown parameter in the linear predictor $\eta
(x,\bftheta)=h(\mathrm{E}(y))$, where $h$ is the link function.

In approximate design context, a design $\xi$ with at most $q$ design
points can be written as $\xi=\{(x_{i},\omega_{i})\}_{i=1}^{q}$, where
$x_i\in[L,U], \omega_i\ge0, i=1,\ldots,q$, $x_i$'s and $\omega_i$'s
are the design points and corresponding design weights, and $\sum_{i=1}^{q} \omega_i=1$. If the weight of a certain design point is
positive, then that design point is a \emph{support point} of the
design, and the number of support points is the \emph{support size} of
the design.

Under the assumption of independent responses, the Fisher information
matrix for $\bftheta$ under design $\xi$ can be written as $ n\sum_{i=1}^q \omega_i \bfM_{x_i}(\bftheta)$, where $n$ is the total sample
size and $\bfM_{x_i}(\bftheta)$ is the information matrix of a single
observation at $x_i$. Since $n$ is only a multiplicative factor, we
prefer using the normalized information matrix, which is $\bfM_{\xi
}(\bftheta)=\sum_{i=1}^q \omega_i \bfM_{x_i}(\bftheta)$.

An optimal design is a design that maximizes the Fisher information
matrix $\bfM_{\xi}(\bftheta)$ under a certain criterion $\Phi$. In this
paper, we focus on a general class of differentiable optimality
criteria. Specifically, let NND($d$) be the set of all \mbox{$d\times d$}
nonnegative definite matrices, PD($d$) be the set of all $d\times d$
positive definite matrices, and $\Phi$ be any function defined on
NND($d$) that satisfies Assumption~\ref{asmp:assumption A} below [see
\citet{p93}, page~115].

\begin{assumption}\label{asmp:assumption A}
Suppose the optimality criterion $\Phi$ is a nonnegative, nonconstant
function defined on NND($d$) such that:
\begin{longlist}[(1)]
\item[(1)] it is concave, that is, $\Phi(\alpha\bfM_1+(1-\alpha)\bfM
_2)\ge\alpha\Phi(\bfM_1)+(1-\alpha)\Phi(\bfM_2)$, where $\alpha\in
(0,1), \bfM_1,\bfM_2\in$ NND($d$);
\item[(2)] it is isotonic, that is, $\Phi(\bfM_1) \ge\Phi(\bfM_2)$ if
$\bfM_1\ge\bfM_2$ under the Loewner ordering, $\bfM_1,\bfM_2\in$ NND($d$);
\item[(3)] it is smooth on PD($d$). By smooth, we mean the function is
differentiable and the first-order partial derivatives are continuous
[for matrix differentiation, $\Phi$~is to be interpreted as a function
of the $d(d+1)/2$-dimensional vector of elements in the upper triangle
of $\bfM$].
\end{longlist}
\end{assumption}

A design $\xi^*$ is $\Phi$-optimal if it maximizes $\Phi(\bfM_{\xi
}(\bftheta))$ with respect to $\xi$.

This class of optimality criteria is very broad and includes, for
example, the well-known $\Phi_p$-optimality criteria with $-\infty<p\le
1$, which are defined as follows. Suppose we are interested in
estimating a smooth function of $\bftheta$, say $g(\bftheta)\dvtx \mathbb
{R}^d \rightarrow\mathbb{R}^v$, where $v \le d$ and $\bfK(\bftheta
)=(\partial g(\bftheta)/\partial\bftheta)^{T}$ has full column rank
$v$. It can be estimated as long as the columns of $\bfK(\bftheta)$ are
contained in the range of $\bfM_{\xi}(\bftheta)$. The information
matrix for $g(\bftheta)$ under design $\xi$ is then defined as $\bfI
_{\xi}(\bftheta)=(\bfK(\bftheta)^{T}\bfM_{\xi}(\bftheta)^{-}\bfK
(\bftheta))^{-1}$, where $\bfM_{\xi}(\bftheta)^{-}$ is a generalized
inverse if $\bfM_{\xi}(\bftheta)$ is singular. Then a $\Phi_p$-optimal
design for $g(\bftheta)$ is defined to maximize
\[
\Phi\bigl(\bfM_{\xi}(\bftheta)\bigr)=\Phi_p\bigl(
\bfI_{\xi}(\bftheta)\bigr)= \biggl(\frac{1}v \operatorname{ trace
}\bigl(\bfI_{\xi}^p(\bftheta)\bigr) \biggr)^{1/p},\qquad
p\in(-\infty,1].
\]
However, $E$-optimality where $g(\bftheta)=\bftheta$ and $p=-\infty$,
is not included here since generally it does not satisfy the smoothness
condition on PD($d$); a short discussion about this can be found in
Section~\ref{sec:Disscusion}. In addition to the $\Phi_p$-optimality
criteria, our general $\Phi$-optimality criteria also include compound
optimality criteria, criteria for evaluating a mixture of information
matrices obtained from nested models [see \citet{p93}, Chapter~11] and so on.

\subsection{Preliminary results}\label{subsec:locally_sub1}
While finding optimal designs is an optimization problem, the
dimensionality of the optimization problem is unknown since the number
of design points, $q$, is unknown. However, it has been observed in the
literature that optimal designs are often saturated designs. This
phenomenon was first discovered in \citet{d54}, and was
generalized to a class of models in
\citeauthor{ys12} (\citeyear{ys09,ys12}), \citet{y10} as well as in \citet{dm11}
and \citet{ds13}, where the latter two papers provided a
different perspective on this phenomenon using \emph{Chebyshev
systems}. Based on these results, optimal designs can be found in a
small \emph{complete class} of designs, denoted as $\Xi$, and in many
cases $\Xi$ only consists of designs with at most $d$ design points.
Here, we briefly introduce a fundamental theorem from \citet{ys12} for our later use. Using the techniques there, we decompose the
Fisher information matrix in the following way (an example is given at
the end of Section~\ref{subsec:locally_sub1}):
%
\begin{equation}
\label{equ:infodecomp} \bfM_{\xi}(\bftheta)=\bfP(\bftheta)\bfC_{\xi}(
\bftheta)\bfP(\bftheta )^T,\qquad \bfC_{\xi}(\bftheta)= \Biggl(\sum
_{i=1}^{q} \omega_i \bfC (
\bftheta,c_i) \Biggr),
\end{equation}
where $\bfC(\bftheta,c)$ is a $d\times d$ symmetric matrix,
\[
\bfC(\bftheta,c)=\pmatrix{
\Psi_{11}(\bftheta,c) & & &
\vspace*{2pt}\cr
\Psi_{21}(\bftheta,c) & \Psi_{22}(\bftheta,c) & &
\vspace*{2pt}\cr
\vdots& \vdots& \ddots&
\vspace*{2pt}\cr
\Psi_{d1}(\bftheta,c) & \Psi_{d2}(\bftheta,c) & \cdots& \Psi
_{dd}(\bftheta,c) },
\]
$\bfP(\bftheta)$ is a $d\times d$ nonsingular matrix that only depends
on $\bftheta$, and $c\in[A,B]$ is a smooth monotonic transformation of
$x$ that depends on $\bftheta$. For the sake of simplicity, we drop
$\bftheta$ from the notation of matrix $\bfC(\bftheta,c)$ and its
elements hereafter [in fact, in many cases a nice decomposition can be
found so that $\bfC(\bftheta,c)$ only depends on $\bftheta$ through
$c$, and $\bftheta$ becomes redundant in the notation].

For some $d_1$, $1\le d_1<d$, define $\bfC_{22}(c)$ as the lower
$d_1\times d_1$ principal submatrix of $\bfC(c)$, that is,
\[
\bfC_{22}(c)=\pmatrix{
\Psi_{d-d_1+1,d-d_1+1}(c) & \cdots & \Psi_{d-d_1+1,d}(c)
\vspace*{2pt}\cr
\vdots& \ddots & \vdots
\vspace*{2pt}\cr
\Psi_{d,d-d_1+1}(c) & \cdots& \Psi_{dd}(c) }.
\]
Choose a maximal set of linearly independent nonconstant functions from
the first $d-d_1$ columns of the matrix $\bfC(c)$, let the number of
functions in this set be $k-1$, and rename them as $\Psi_\ell(c)$, $\ell
=1,\ldots,k-1$. Let $\Psi_k(c)=\bfC_{22}(c)$, and define the functions
$f_{\ell,t}(c)$, $1 \le t \le\ell\le k$, to be

{\fontsize{9}{10}{\selectfont{
\[
\pmatrix{ f_{1,1}=\Psi_1'
& & & & & \vspace*{2pt}
\cr
f_{2,1}=\Psi_2' &\displaystyle
f_{2,2}= \biggl(\frac{f_{2,1}}{f_{1,1}} \biggr)' & & & &
\vspace*{2pt}
\cr
f_{3,1}=\Psi_3' &
\displaystyle f_{3,2}= \biggl(\frac{f_{3,1}}{f_{1,1}} \biggr)' &
\displaystyle f_{3,3}= \biggl(\frac{f_{3,2}}{f_{2,2}} \biggr)' & & &
\vspace*{2pt}
\cr
f_{4,1}=\Psi_4' &
\displaystyle f_{4,2}= \biggl(\frac{f_{4,1}}{f_{1,1}} \biggr)' &
\displaystyle f_{4,3}= \biggl(\frac{f_{4,2}}{f_{2,2}} \biggr)' &
\displaystyle f_{4,4}= \biggl(\frac
{f_{4,3}}{f_{3,3}} \biggr)' & &
\vspace*{2pt}
\cr
\vdots& \vdots & \vdots & \vdots & \ddots& \vspace*{2pt}
\cr
f_{k,1}=\Psi_k' & f_{k,2} = \biggl(
\displaystyle \frac{f_{k,1}}{f_{1,1}} \biggr)' & f_{k,3}= \biggl(
\displaystyle \frac{f_{k,2}}{f_{2,2}} \biggr)' & f_{k,4}= \biggl(
\frac
{f_{k,3}}{f_{3,3}} \biggr)' & \hspace*{-1pt}\cdots\hspace*{-1pt}& f_{k,k}= \biggl(
\displaystyle \frac
{f_{k,k-1}}{f_{k-1,k-1}} \biggr)' }\hspace*{-2pt}
, \label{df}
\]}}}
\hspace*{-5pt}where the entries in the last row are matrices, and the derivatives of
matrices are element-wise derivatives (assuming all derivatives exist).
Define matrix $\bfF(c) = \prod_{\ell=1}^k f_{\ell,\ell}(c)$. Then the
following theorem due to \citet{ys12} is available [see also
\citet{ds13}, Theorem~3.1].

\begin{thmm}[{[\citet{ys12}]}]\label{thmm:thmm1}
For a regression model with a single covariate, suppose that either
$\bfF(c)$ or $-\bfF(c)$ is positive definite for all $c \in[A,B]$.
Then the following results hold:
\begin{longlist}[(a)]
\item[(a)] If $k=2m-1$ is odd and $\bfF(c)<0$, then designs with at
most $m$ design points, including point $A$, form a complete class $\Xi$.
\item[(b)] If $k=2m-1$ is odd and $\bfF(c)>0$, then designs with at
most $m$ design points, including point $B$, form a complete class $\Xi$.
\item[(c)] If $k=2m$ is even and $\bfF(c)<0$, then designs with at most
$m$ design points, form a complete class $\Xi$.
\item[(d)] If $k=2m-2$ is even and $\bfF(c)>0$, then designs with at
most $m$ design points, including both $A$ and $B$, form a complete
class $\Xi$.
\end{longlist}
\end{thmm}

It is helpful to sketch how Theorem~\ref{thmm:thmm1} is proved. For
some carefully chosen $d_1$ (see example below) where one of the
conditions in Theorem~\ref{thmm:thmm1} holds, it can be proved that for
any design $\xi\notin\Xi$, we can find a design $\tilde{\xi}\in\Xi$
such that $\bfC_{\tilde{\xi}}(\bftheta)\ge\bfC_{\xi}(\bftheta)$ under
the Loewner ordering, hence $\bfM_{\tilde{\xi}}(\bftheta)\ge\bfM_{\xi
}(\bftheta)$. To be specific, $\bfC_{\tilde{\xi}}(\bftheta)-\bfC_{\xi
}(\bftheta)$ has a \emph{positive definite} lower $d_1\times d_1$
principal submatrix, and is 0 everywhere else. So the search for
optimal designs can be restricted within $\Xi$.

Theorem~\ref{thmm:thmm1} also applies to generalized linear models.
Besides, while it is stated in terms of the ``transformed design
point'' $c$, the result can be easily translated back into $x$ using
the relationship between them, and we will state results in $x$ unless
otherwise specified.

In Theorem~\ref{thmm:thmm1}, there are four different types of complete
classes, the difference being whether one or both of the endpoints are
fixed design points (note however a fixed design point can have weight
0 so that it need not be a support point). To make it easier to
distinguish, let $\operatorname{fix}(\Xi)$ denote the set of fixed design points for
the designs in the complete class $\Xi$. For example, $\operatorname{fix}(\Xi
)=\varnothing$ and $\{L,U\}$ refers to the complete classes in
Theorem~\ref{thmm:thmm1}(c) and (d), respectively.

Applications of Theorem~\ref{thmm:thmm1} can be found in
\citeauthor{ys12} (\citeyear{ys09,ys12}) and \citet{y10}. Obviously,
$m\ge d$, however, in many applications we actually find $m=d$. Take
the LINEXP model from \citet{ys12} as an example.

The LINEXP model is used to characterize tumor growth delay and
regrowth. The natural logarithm of tumor volume is modeled using a
nonlinear regression model with mean
%
\begin{equation}
\label{equ:modlinexp} \eta(x,\bftheta)=\theta_1+\theta_2e^{\theta_3x}+
\theta_4x,
\end{equation}
where $x\in[L, U]$ is the time, $\theta_1+\theta_2$ is the logarithm
of initial tumor volume, $\theta_3<0$ is the rate at which killed cells
are eliminated, $\theta_4>0$ is the final growth rate.

The information matrix for $\bftheta$ can be written in the form of
(\ref{equ:infodecomp}) with
\[
P(\bftheta)=\pmatrix{ 1 & 0 & 0 & 0
\vspace*{2pt}\cr
0 & 1 & 0 & 0
\vspace*{2pt}\cr
0 & 0 & 0 & {\theta_2}/{\theta_3}
\vspace*{2pt}\cr
0 & 0 & {1}/{\theta_3} & 0 },\qquad \bfC(c)=
\pmatrix{1 & & &
\vspace*{2pt}\cr
e^c & e^{2c} & &
\vspace*{2pt}\cr
c & ce^c & c^2 &
\vspace*{2pt}\cr
ce^c & ce^{2c} & c^2e^c &
c^2e^{2c}},
\]
where $c=\theta_3x\in[A,B]=[\theta_3U,\theta_3L]$. Let $d_1=2$,
$\textbf{C}_{22}(c)$ be the lower $2\times2$ principal submatrix of
$\textbf{C}(c)$, and $\Psi_1(c)=c,\Psi_2(c)=e^c,\Psi_3(c)= ce^c, \Psi
_4(c)=e^{2c}, \Psi_5(c)= ce^{2c}$ be the set of linearly independent
nonconstant functions from the first two columns of $\textbf{C}(c)$.
Then $k=6$, $f_{1,1}=1, f_{2,2}=e^c, f_{3,3}=1, f_{4,4}=4e^c,
f_{5,5}=1$, and
\[
f_{6,6}(c)= %
\pmatrix{ 2e^{-2c} & e^{-c}/2
\vspace*{2pt}
\cr
e^{-c}/2 & 2 } %
,\qquad \bfF(c)=\prod
_{\ell=1}^6 f_{\ell,\ell}(c)= %
\pmatrix{
8 & 2e^{c}\vspace*{2pt}
\cr
2e^{c} & 8e^{2c} }.
\]
Because $\bfF(c)>0$, Theorem~\ref{thmm:thmm1}(d) can be applied with
$m=4=d$, and $\Xi$ consists of designs with at most four design points
including both endpoints, thus $\operatorname{fix}(\Xi)=\{L,U\}$.

\subsection{Identifying the optimal design}\label{subsec:locally_sub2}
If one of the cases in Theorem~\ref{thmm:thmm1} holds, an optimal
design exists of the form $\xi=\{(x_i, \omega_i)\}_{i=1}^{m}$, where
$x_i$'s are strictly increasing, with $x_1$ or $x_m$ possibly fixed to
be $L$ or $U$, respectively; $\omega_i$'s are nonnegative, and $\omega
_1=1-\sum_{i=2}^m \omega_i$. Let $\bfZ$ be the vector of unknown design
points (i.e., exclude $x_1$ or $x_m$ if fixed to be the endpoint) and
$m-1$ unknown weights $\omega_2,\ldots,\omega_m$. For example, for the
LINEXP model in (\ref{equ:modlinexp}), $\bfZ=(x_2, x_3, \omega_2,\omega
_3, \omega_4)^T$ since $m=4$ and $x_1=L, x_4=U$. Thus we can use $\bfZ$
to represent the design $\xi$. Now the objective function $\Phi(\bfM
_{\xi}(\bftheta))$ is a function of $\bfZ$, denoted as $\tilde{\Phi
}(\bfZ)$, and it is smooth by the smoothness of $\Phi$. To find an
optimal design, we need to maximize $\tilde{\Phi}(\bfZ)$ with respect
to $\bfZ$. The simplest way is to find the critical points,
specifically, the \emph{feasible critical points}, as defined below.

\begin{definition}\label{def:feas_critic}
A critical point of $\tilde{\Phi}(\bfZ)$, $\bfZ^c$, is a feasible
critical point if all the design points in $\bfZ^c$ are within $[L,U]$
and all $m-1$ weights are positive with summation less than 1.
\end{definition}

With Definition $\ref{def:feas_critic}$, a feasible critical point
gives a design with $m$ support points. Moreover, Theorem~\ref
{thmm:thmm2} states the conditions such that a feasible critical point
gives a globally optimal design.

\begin{thmm}\label{thmm:thmm2}
Assume one of the cases in Theorem~\ref{thmm:thmm1} holds, then for any
feasible critical point of $\tilde{\Phi}(\bfZ)$, its corresponding
design is a $\Phi$-optimal design.
\end{thmm}

\begin{pf}
See the \hyperref[sec:Proofs]{Appendix}.
\end{pf}

Theorem~\ref{thmm:thmm2} gives an implicit solution of an optimal
design if there exists a feasible critical point. Such a point can be
given explicitly in special situations, but not in general due to the
complexity of the objective function. Nevertheless, we have an implicit
solution and it can be easily solved using Newton's algorithm. However,
we need to guarantee the existence of a feasible critical point in the
first place. Theorem~\ref{thmm:thmm3} gives some sufficient conditions
that a feasible critical point exists.

\begin{thmm}\label{thmm:thmm3}
Suppose one of the cases in Theorem~\ref{thmm:thmm1} holds and any $\Phi
$-optimal design has at least $m$ support points. Further assume one of
the following four conditions holds:
\begin{longlist}[(a)]
\item[(a)] $\operatorname{fix}(\Xi)=\{L\}$, and the information matrix $\bfM
_{U}(\bftheta)$ is 0;
\item[(b)] $\operatorname{fix}(\Xi)=\{U\}$, and $\bfM_{L}(\bftheta)=0$;
\item[(c)] $\operatorname{fix}(\Xi)=\varnothing$, and $\bfM_{U}(\bftheta)=\bfM
_{L}(\bftheta)=0$;
\item[(d)] $\operatorname{fix}(\Xi)=\{L,U\}$.
\end{longlist}
Then a feasible critical point of $\tilde{\Phi}(\bfZ)$ must exist, and
by Theorem~\ref{thmm:thmm2}, any such point gives a $\Phi$-optimal design.
\end{thmm}

\begin{pf}
Let $\xi^*\in\Xi$ be a $\Phi$-optimal design, then $\xi^*$ has at
least $m$ support points. By Theorem~\ref{thmm:thmm1}, designs in the
complete class have at most $m$ support points, hence $\xi^*$ has
exactly $m$ support points. Let $\bfZ^*$ be the vector corresponding to
$\xi^*$ according to the definition in the beginning of Section~\ref{subsec:locally_sub2}. For each of conditions (a)$\sim$(d), we know the
design points in $\bfZ^*$ do not include any of the endpoints (recall
the fixed design points are excluded in $\bfZ^*$), hence they all
belong to the open interval $(L,U)$. The weights in $\bfZ^*$ are all
positive, hence all belong to the open interval $(0,1)$, so $\bfZ^*$ is
not on the boundary and must be a critical point of $\tilde{\Phi}(\bfZ
)$. This proves the existence.
\end{pf}

The condition in Theorem \ref{thmm:thmm3} that
every $\Phi$-optimal design has at least $m$ support points is met with
$m=d$ for many models and optimality criteria. For example, when
$K(\theta)$ is a nonsingular matrix, any $\Phi$-optimal design has at
least $d$ support points for commonly used  optimality criteria. On the
other hand, as we have stated, for many models, the complete class given
by Theorem \ref{thmm:thmm1} only consists of designs with at most $d$
support points.
The
condition (d) is found to be satisfied for several models, as we will
see in Section~\ref{sec:app}. For condition (a), usually $\bfM
_{U}(\bftheta)=0$ only when $U=\infty$, so the condition fails if we
are interested in a finite design region, and so do conditions~(b) and
(c). This issue will be addressed later in Theorem~\ref{thmm:thmm5}.

Useful results can be obtained by applying Theorem~\ref{thmm:thmm3} to
the most commonly used $\Phi_p$-optimality criteria. In particular, we
are interested in $\Phi_p$-optimal designs for $\bftheta$ or $\bfa
^T\bftheta$, where $\bfa=(a_1,\ldots,a_d)^T$ is a $d\times1$ vector
such that $\bfa^T\bftheta$ is only estimable with at least $d$ support
points. Adopting the notation in \citet{kw65}, define
\[
A^*=\bigl\{\bfa| \bfa^T\bftheta\mbox{ is only estimable with at
least } d \mbox{ support points} \bigr\}.
\]
Now Corollary~\ref{corl:corl1} gives applications of Theorem~\ref
{thmm:thmm3} to $\Phi_p$-optimal designs.

\begin{corollary}\label{corl:corl1}
Suppose that one of the cases in Theorem~\ref{thmm:thmm1} holds with
$m=d$, and one of the four conditions in Theorem~\ref{thmm:thmm3} is
met. Consider $\Phi_p$-optimal design for $g(\bftheta)$ where
$g(\bftheta)$ satisfies either case \textup{(i)} or \textup{(ii)} below:
\begin{longlist}[(ii)]
\item[(i)] $g(\bftheta)=\bftheta$ or a reparameterization of $\bftheta$;
\item[(ii)] $g(\bftheta)=\bfa^T\bftheta, \bfa\in A^*$.
\end{longlist}
Then a feasible critical point of $\tilde{\Phi}(\bfZ)$ exists, and any
such point gives a $\Phi_p$-optimal design for $g(\bftheta)$.
\end{corollary}

\begin{remark}\label{rmk:rmk1}
In Corollary~\ref{corl:corl1}(i), a special case of a
reparameterization is $g(\bftheta)=\bfW\bftheta$, where $\bfW$ is a
diagonal matrix with positive diagonal elements. This makes cov$(g(\hat
{\bftheta}))$ a rescaled version of cov$(\hat{\bftheta})$, and it makes
sense when $\operatorname{var}(\hat{\theta}_i)$'s are of different orders of
magnitude. For example, in \citet{d97}, the author proposed
``standardized'' optimality criteria, where the matrix $\bfW$ has
diagonal elements $\bfW_{\mathit{ii}}=\sqrt{1/(\bfM^{-1}_{\xi^*_i})_{\mathit{ii}} }$, $\xi
_i^*$ is the $c$-optimal design for estimating $\theta_i$ alone,
$i=1,\ldots,d$. Under the conditions of Corollary~\ref{corl:corl1}(i),
finding such optimal designs is easy after we find $\xi^*_i$'s.
\end{remark}

\begin{remark}\label{rmk:rmk2}
Corollary~\ref{corl:corl1}(ii) considers $c$-optimality. When $\bfa\in
A^*$, the $c$-optimal design is supported at the full set of Chebyshev
points in many cases [see \citet{s68}], but our method gives another
way of finding $c$-optimal designs. When $\bfa\notin A^*$, sometimes a
feasible critical point still exists, and it still gives an optimal
design. However, if there is no such critical point, then the
$c$-optimal design must be supported at fewer points, which may not be
the Chebyshev points, and this problem becomes harder. Nevertheless, we
can approximate such $c$-optimal designs. Suppose $a_1 \neq0$,
consider $g_{\epsilon}(\bftheta)=(\bfa^T\bftheta, \epsilon\theta_2,
\ldots, \epsilon\theta_d)^T, \epsilon>0$. A $\Phi_p$-optimal design
for $g_\epsilon(\bftheta)$ can be found easily by Corollary~\ref
{corl:corl1}(i). Let $\epsilon\rightarrow0$, it can be shown that
these $\Phi_{p}$-optimal designs will eventually converge to the
$c$-optimal design for $\bfa^T\bftheta$ (i.e., the efficiencies of
these $\Phi_{p}$-optimal designs under $c$-optimality will converge to
1), for any $p\le-1$. Some examples are provided in Section~\ref{subsec:app_sub2}.
\end{remark}

To verify the condition $\bfa\in A^*$, let $\bff(x,\bftheta
)=(f_1(x,\bftheta), \ldots, f_{d}(x,\bftheta))= \break \partial\eta(x,\bftheta
)/\partial\bftheta$. The condition $\bfa\in A^*$ is equivalent to
%
\begin{equation}
\label{equ:cheby} %
\left| \matrix{ f_1(x_1,\bftheta)
& \cdots& f_1(x_{d-1},\bftheta) & a_1
\vspace*{2pt}
\cr
f_2(x_1,\bftheta) & \cdots&
f_2(x_{d-1},\bftheta) & a_2 \vspace*{2pt}
\cr
\vdots& \ddots & \vdots& \vdots \vspace*{2pt}
\cr
f_d(x_1,
\bftheta) & \cdots& f_d(x_{d-1},\bftheta) & a_d
}\right| %
\neq0
\end{equation}
for all $L\le x_1<x_2<\cdots<x_{d-1}\le U$ (this is also true for
generalized linear models). In particular, if we are interested in
estimating the individual parameter $\theta_i$, that is, $\bfa=\bfe_i$
where $\bfe_i=(0,\ldots,0,1,0,\ldots,0)^T$ denotes the $i$th unit
vector, then $\bfe_i \in A^*$ is equivalent to $\bff_{-i}=\{f_j| j\in\{
1,\ldots,d\} \setminus\{i\} \}$ being a Chebyshev system [see \citet{ks66}], which is easier to verify. Here, the traditional
definition of a Chebyshev system is used, which only requires the
determinant in~(\ref{equ:cheby}) to be \emph{nonzero} instead of positive.

Next, the uniqueness of optimal designs can also be established under
mild conditions. We first introduce some additional terminology. A
criterion $\Phi$ is called \emph{strictly isotonic} on PD($d$) if
\[
\Phi(\bfM_1)>\Phi(\bfM_2)\qquad \mbox{for any }
\bfM_1\ge\bfM_2>0 \mbox{ and } \bfM_1\neq
\bfM_2.
\]
It is called \emph{strictly concave} on PD($d$) if
\begin{eqnarray}
\Phi\bigl(\alpha\bfM_1+(1-\alpha)\bfM_2\bigr)>\alpha
\Phi(\bfM_1)+(1-\alpha)\Phi (\bfM_2),\nonumber
\\
\eqntext{\mbox{for any } \alpha\in(0,1), \bfM_1>0, \bfM_2\ge0
\mbox{ and } \bfM_2\not\propto\bfM_1.}
\end{eqnarray}
For example, $\Phi_p$-optimality criteria are both strictly isotonic
and strictly concave on PD($d$) when $g(\bftheta)$ is $\bftheta$ or a
reparameterization of $\bftheta$ and $p \in(-\infty,1)$ [see
\citet{p93}, page~151]. Moreover, a compound optimality criterion
which involves a strictly isotonic and strictly concave criterion is
also strictly isotonic and strictly concave. For these criteria, we
have Theorem~\ref{thmm:thmm4}.

\begin{thmm}\label{thmm:thmm4}
Assume that one of the cases in Theorem~\ref{thmm:thmm1} holds. If $\Phi
$ is both strictly isotonic and strictly concave on PD($d$) and there
exists a $\Phi$-optimal design $\xi^*$ which has at least $d$ support
points, then $\xi^*$ is the unique $\Phi$-optimal design. In
particular, the $\Phi_p$-optimal design under Corollary~\ref
{corl:corl1}\textup{(i)} is unique for $p\in(-\infty,1)$.
\end{thmm}

\begin{pf}
See the \hyperref[sec:Proofs]{Appendix}.
\end{pf}

\begin{remark}\label{rmk:rmk3}
The $c$-optimality criterion with $g(\bftheta)=\bfa^T\bftheta$ maybe
neither strictly concave nor strictly isotonic on PD($d$). However, if
$\bfa\in A^*$ and $\bff(x,\bftheta)$ is a Chebyshev system, the
uniqueness is proved in \citet{s68}.
\end{remark}

The uniqueness is not only of interest in itself, but also has
implications for finding optimal designs. As we have stated earlier,
conditions~(a), (b) and (c) in Theorem~\ref{thmm:thmm3} may only hold
on a large design region, call it the \emph{full design region}. Let
$\xi^{**}$ be a $\Phi$-optimal design on the full design region with
smallest support point $x_{\mathrm{min}}^{**}$ and largest support point
$x_{\mathrm{max}}^{**}$. Then for a smaller design region $[L,U]$, under the
same optimality criterion $\Phi$, we have Theorem~\ref{thmm:thmm5}.

\begin{thmm}\label{thmm:thmm5}
Assume that one of the cases in Theorem~\ref{thmm:thmm1} holds for the
full design region, and both $\Phi$-optimal designs on $[L,U]$ and the
full design region are unique with support size $m$, then we have:
\begin{longlist}[(a)]
\item[(a)] under $\operatorname{fix}(\Xi)=\{L\}$, if $U<x_{\mathrm{max}}^{**}$, then the $\Phi
$-optimal design on $[L,U]$ has both $L$ and $U$ as support points;
otherwise, the optimal design is $\xi^{**}$;
\item[(b)] under $\operatorname{fix}(\Xi)=\{U\}$, if $x_{\mathrm{min}}^{**}<L$, then the $\Phi
$-optimal design on $[L,U]$ has both $L$ and $U$ as support points;
otherwise, the optimal design is $\xi^{**}$;
\item[(c)] under $\operatorname{fix}(\Xi)=\varnothing$, if $x_{\mathrm{min}}^{**}<L$ or
$U<x_{\mathrm{max}}^{**}$, then the $\Phi$-optimal design on $[L,U]$ has at
least one endpoint as a support point; otherwise, the optimal design is
$\xi^{**}$.
\end{longlist}
\end{thmm}

\begin{pf}
We only give the proof for case (a), others being similar. When $U\ge
x_{\mathrm{max}}^{**}$, the design $\xi^{**}$ is still a feasible design on the
region $[L,U]$, and it is optimal because it is optimal on the full
design region. When $U<x_{\mathrm{max}}^{**}$, $\xi^{**}$ is no longer a
feasible design, let $\xi^*$ be the optimal design on $[L,U]$. A
complete class of the same type exists for design region $[L,U]$
because, for example, $\bfF(c)>0$ on the full design region implies
$\bfF(c)>0$ on the smaller design region. So $x_1^*=L$. If the largest
support point $x_{m}^*<U$, then $\bfZ^*=(x_2^*,\ldots,x_{m}^*, \omega
_2^*,\ldots,\omega_m^*)^T$ must be a critical point of $\tilde{\Phi
}(\bfZ)$. Now if we consider the optimal design problem on the full
design region again, $\bfZ^*$ is a feasible critical point, and by
Theorem~\ref{thmm:thmm2}, $\xi^*$ must be an optimal design on the full
design region. However, $\xi^*\neq\xi^{**}$, this contradicts the
uniqueness assumption.
\end{pf}

\section{Application}\label{sec:app}
The theorems we have established can be used to find optimal designs
for many models. In Sections~\ref{subsec:app_sub1} through \ref
{subsec:app_sub3}, we consider $\Phi_p$-optimal designs for models with
two, three and four or six parameters, respectively. In Section~\ref{subsec:app_sub4}, we consider polynomial regression models with
arbitrary $d$ parameters under more general optimality criteria.

\subsection{Models with two parameters}\label{subsec:app_sub1}
Yang and Stufken (\citeyear{ys09}) considered complete class results for
two-parameter models, including logistic/probit regression model,
Poisson regression model and Michaelis--Menten model. The theorems we
have established can be used to find the optimal designs. Take the
Poisson regression model as an example (the applications to other
models are similar). It has the following form:
\[
\eta(x,\bftheta)=\log\bigl(\mathrm{E}(y)\bigr)=\theta_1+
\theta_2 x,\qquad x\in[L,U].
\]
Theorem~\ref{thmm:thmm1}(b) can be applied to this model, and a
complete class consists of designs with at most 2 design points
including one boundary point [see \citet{ys09}, Theorem~4].
Specifically, when $\theta_2>0$, $U$ is a fixed design point, and $\bfM
_{-\infty}(\bftheta)=0$ [since $\bfM_{x}(\bftheta)=e^{\theta_1+\theta
_2x}(1,x)^T(1,x)$]; when $\theta_2<0$, $L$ is a fixed design point, and
$\bfM_{\infty}(\bftheta)=0$. Thus, on any one-sided restricted region
$(-\infty, U]$ (when $\theta_2>0$) or $[L,\infty)$ (when $\theta_2<0$),
$\Phi_p$-optimal designs for $\bftheta$ can be found by solving for the
critical points, according to Corollary~\ref{corl:corl1}(i). For
$c$-optimality, recall $\bff(x,\bftheta)=\partial\eta(x,\bftheta
)/\partial\bftheta=(1,x)$, thus $\bff_{-2}=\{1\}$ is a Chebyshev
system, which means $\theta_2$ can only be estimated with at least
$d=2$ support points. Therefore, according to Corollary~\ref
{corl:corl1}(ii), an $\bfe_2$-optimal design ($c$-optimal design for
$\theta_2$) can also be found by solving for the critical points.

In particular, $D$- and $\bfe_2$-optimal designs can be found
analytically through symbolic computation software (e.g., by using the
solve function in Matlab) and are listed in~(\ref{equ:doptpois}) and
(\ref{equ:e2optpois}). Note that they do not depend on $\theta_1$ since
$e^{\theta_1}$ is merely a multiplicative factor in $\bfM_x(\bftheta)$:
%
\begin{eqnarray}
\label{equ:doptpois} \xi_D^*&=&\cases{ %
\bigl
\{(U-2/\theta_2,1/2), (U,1/2)\bigr\}, & \quad$\theta_2>0,$
\vspace*{2pt}\cr
\bigl\{(L-2/\theta_2,1/2), (L,1/2)\bigr\}, &\quad $\theta_2<0,$}
\\
\label{equ:e2optpois} \xi_{\bfe_2}^*&=&\cases{ %
\bigl
\{(U-2.557/\theta_2,0.782), (U,0.218)\bigr\}, &\quad  $\theta_2>0,$
\vspace*{2pt}\cr
\bigl\{(L-2.557/\theta_2,0.782), (L,0.218)\bigr\}, &\quad
$\theta_2<0$.}
\end{eqnarray}
However, $A$-optimal designs do not have explicit forms. Nevertheless,
the solutions can be found easily using Newton's algorithm. For the
case of $\theta_2<0$, some examples are listed in Table~\ref{Tab:tab1}
(again the optimal designs do not depend on $\theta_1$).

\begin{table}
\tablewidth=220pt
\caption{$A$-optimal designs for Poisson regression model on $[0,\infty)$}
\label{Tab:tab1}
\begin{tabular*}{220pt}{@{\extracolsep{\fill}}lcc@{}}
\hline
& \multicolumn{2}{c@{}}{${A}${-optimal}} \\[-6pt]
& \multicolumn{2}{c@{}}{\hrulefill} \\
\multicolumn{1}{@{}l}{${\theta_2}$} &
\multicolumn{1}{c}{${(x_1, x_2)}$} &
\multicolumn{1}{c@{}}{${(\omega_1, \omega_2)}$} \\
\hline
$-1$ & (0, 2.261) & (0.444, 0.556) \\
$-2$ & (0, 1.193) & (0.320, 0.680) \\
\hline
\end{tabular*}
\end{table}

In addition, the $\Phi_p$-optimal design for $\bftheta$ and $\bfe
_2$-optimal design are unique, due to Theorem~\ref{thmm:thmm4}. For
finite design regions, Theorem~\ref{thmm:thmm5} can be applied. For
example, the $A$-optimal design for $\bftheta=(1,-1)^T$ on $[0,U]$ when $U
\ge2.261$ is $\{(0,0.444),(2.261,0.556)\}$; when $U<2.261$, the
optimal design is supported at exactly two points $0$ and $U$, and the
weights can be determined easily.

\subsection{Models with three parameters}\label{subsec:app_sub2}
Dette et~al. (\citeyear{dbpp08,dkbb10}) considered optimal designs
for the Emax and log-linear models. These models, often used to model
dose-response curves, are nonlinear regression models with means
\[
\eta(x,\bftheta)=\cases{ %
\theta_1+
\theta_2x/(x+\theta_3), &\quad$ \mathrm{Emax},$
\vspace*{2pt}\cr
\theta_1+\theta_2\log(x+\theta_3), &\quad
$\mathrm{log\mbox{-}linear}.$ }
\]
Here, $x\in[L,U]\subseteq(0,\infty)$ is the dose range, $\theta_2>0$
and $\theta_3>0$. Theorem~\ref{thmm:thmm1}(d) can be applied to both
models, and a complete class consists of designs with at most 3 design
points including \emph{both} endpoints [\citet{y10}, Theorem~3]. Hence,
Corollary~\ref{corl:corl1} is applicable on design space $[L,U]$. In
particular, $D$-optimal designs can be computed explicitly using
symbolic computation software, and are listed in~(\ref{equ:doptEmax}).
They are consistent with the results in \citet{dkbb10}:
%
\begin{equation}
\label{equ:doptEmax} \xi_D^*=\cases{ %
 \bigl
\{(L,1/3),\bigl(x^*_{E},1/3\bigr),(U,1/3)\bigr\}, &\quad $\mathrm{Emax},$
\vspace*{2pt}\cr
\bigl\{(L,1/3),\bigl(x^*_{l},1/3 \bigr), (U,1/3) \bigr\}, &\quad
$\mathrm{log\mbox{-}linear},$}
\end{equation}
where
%
\begin{eqnarray}
\label{equ:xexl} x^*_{E}&=&\frac{L(U+\theta_3)+U(L+\theta_3)}{L+U+2\theta_3},
\nonumber
\\[-8pt]
\\[-8pt]
\nonumber
 x^*_{l}&=&
\frac{(L+\theta_3)(U+\theta_3)}{U-L}\log \biggl(\frac{U+\theta
_3}{L+\theta_3} \biggr)-\theta_3.
\end{eqnarray}

For $A$-optimality, numerical solutions can be obtained easily by
Newton's algorithm. Table~\ref{Tab:tab2} gives some examples for the
Emax model using parameter settings in \citet{dbpp08} (the optimal
designs do not depend on $\theta_1$ since it is not involved in the
information matrix; and although it seems that the optimal weights are
constant, they do change gradually with $\theta_2$ and $\theta_3$).

%
\begin{table}
\caption{$A$-optimal designs for the Emax model on $[0,150]$}\label{Tab:tab2}
\begin{tabular*}{\textwidth}{@{\extracolsep{\fill}}lccc@{}}
\hline
\multicolumn{1}{@{}l}{${\theta_2}$} &
\multicolumn{1}{c}{${\theta_3}$} &
\multicolumn{1}{c}{${(x_1, x_2, x_3)}$} &
\multicolumn{1}{c@{}}{${(\omega_1, \omega_2,
\omega_3)}$} \\
\hline
$7/15$ & 15 & (0, 12.50, 150) & (0.250, 0.500, 0.250)\\
$7/15$ & 25 & (0, 18.75, 150) & (0.250, 0.500, 0.250)\\
$10/15$ & 25 & (0, 18.75, 150) & (0.250, 0.500, 0.250)\\
\hline
\end{tabular*}
\end{table}

For $c$-optimality, \citet{dkbb10} gave explicit solutions for
$ED_p$-optimal designs, where an $ED_p$-optimal design is a design that
is optimal for estimating the dose that achieves $100p\%$ of the
maximum effect in dose range $[L,U]$, $0<p<1$. In fact,
$ED_p$-optimality is equivalent to $\bfe_3$-optimality regardless of
$p$, and we can find the optimal designs using our method. First, we have
\[
\bff(x,\bftheta)=\cases{ %
 \bigl(1,x/(x+
\theta_3), -\theta_2x/(x+\theta_3)^2
\bigr), &\quad $\mathrm{Emax},$
\vspace*{2pt}\cr
\bigl(1,\log(x+\theta_3), \theta_2/(x+
\theta_3)\bigr), &\quad$\mathrm{log\mbox{-}linear}.$}
\]
It is easy to prove for both the Emax and log-linear models that $\bff
_{-3}$ is a Chebyshev system, which means that $\theta_3$ is only
estimable with at least $d=3$ support points. So $\bfe_3$-optimal
designs can be found by solving for the critical points, by
Corollary~\ref{corl:corl1}(ii). The solutions can be found explicitly
using symbolic computation software and are listed in~(\ref
{equ:e3optEmax}). They are consistent with the results in \citet{dkbb10}:
%
\begin{equation}
\label{equ:e3optEmax} \xi_{\bfe_3}^*=\xi^*_{ED_p}=\cases{ %
 \bigl\{(L,1/4),\bigl(x^*_{E},1/2\bigr),(U,1/4)\bigr
\}, &\quad $\mathrm{Emax},$
\vspace*{2pt}\cr
\bigl\{\bigl(L,\omega^*_{l}\bigr), \bigl(x^*_{l}, 1/2
\bigr), \bigl(U,1/2-\omega^*_{l}\bigr) \bigr\}, &\quad$\mathrm{log\mbox{-}linear},$}
\end{equation}
where $x^*_{E}$ and $x^*_{l}$ are the same as in~(\ref{equ:xexl}), and
\[
\omega^*_{l}=\frac{\log(x^*_{l}+\theta_3)-\log(U+\theta_3)}{2(\log
(L+\theta_3)-\log(U+\theta_3))}.
\]

Regarding $\bff_{-2}$, it can be shown that it is always a Chebyshev
system for the log-linear model, and it is a Chebyshev system for the
Emax model if $\theta_3\notin(L,U)$. In such cases, $\bfe_2$-optimal
designs can be found according to Corollary~\ref{corl:corl1}(ii), and
the solutions can be derived analytically as shown in~(\ref{equ:e2optEmax}):
%
\begin{equation}\qquad
\label{equ:e2optEmax} \xi_{\bfe_2}^*=\cases{ %
\displaystyle \biggl\{
\biggl(L,\frac{1}4-\frac{(U-L)\theta_3}{8(\theta_3^2-LU)}\biggr),\biggl(x^*_{E},
\frac{1}2\biggr),\biggl(U,\frac{1}4+\frac{(U-L)\theta_3}{8(\theta_3^2-LU)}\biggr)
\biggr\}, \vspace*{2pt}\cr
\qquad \hspace*{185pt}\mathrm{Emax}, \theta_3\notin(L,U),
\vspace*{8pt}\cr
\displaystyle\biggl\{\biggl(L,\frac{(U-x^*_{l})(L+\theta_2)}{2(U-L)(x^*_{l}+\theta_2)}\biggr), \biggl(x^*_{l},
\frac{1}2\biggr), \biggl(U,\frac{(x^*_{l}-L)(U+\theta
_2)}{2(U-L)(x^*_{l}+\theta_2)}\biggr) \biggr\},\vspace*{2pt}\cr
\qquad\hspace*{224pt}
\mathrm{log\mbox{-}linear}.}\hspace*{-12pt}
\end{equation}

When $\theta_3\in(L,U)$, $\bff_{-2}$ is no longer a Chebyshev system
for the Emax model. However, if $|(U-L)\theta_3|<|2(\theta_3^2-LU)|$,
the weights of $\xi_{\bfe_2}^*$ in~(\ref{equ:e2optEmax}) are still
positive, and the design is still $\bfe_2$-optimal; otherwise, the
optimal design is supported at fewer than 3 points, which may not be
the Chebyshev points. Nevertheless, we can approach the optimal design
using the method in Remark~\ref{rmk:rmk2}. To show this, consider the
setting where the dose range is $[0,150]$, $\theta_2=7/15$ and $\theta
_3=25$. The exact $\bfe_2$-optimal design can be found to be $\xi_{\bfe
_2}^*=\{(\theta_3^2/U, 0.5),(U,0.5)\}=\{(25/6,0.5),(150,0.5)\}$ using
Elfving's method [\citet{e52}]. Now let $\epsilon=10^{-5},10^{-6},10^{-7}$; the $\Phi
_p$-optimal designs for estimating $g_{\epsilon}(\bftheta)=(\epsilon
\theta_1, \theta_2,\epsilon\theta_3)^T$ can be found by Corollary~\ref
{corl:corl1}(i) and are used to approximate the $\bfe_2$-optimal
design. Table~\ref{Tab:tab3} shows the errors and $1-$ efficiencies of
the approximation for $p=-1$ and $-3$. As we can see, the error gets
sufficiently small after a few iterations, especially when $|p|$ is
larger; however, due to singularity issues, the error cannot be made
arbitrary small.

\begin{table}
\tabcolsep=0pt
\caption{Approximating two point $\bfe_2$-optimal design using three
point designs for the Emax model}\label{Tab:tab3}
\begin{tabular*}{\textwidth}{@{\extracolsep{\fill}}lcccccc@{}}
\hline
${p}$ & ${\epsilon}$ & ${|x_2-\frac{25}{6}|/\frac{25}{6}}$ &
${|\omega_1|}$ &
${|\omega_2-0.5|}$ & ${|\omega_3-0.5|}$ & ${1-\operatorname{eff}}$ \\
\hline
$-1$ 
& $10^{-5}$ & $10^{-2}$ & $4\cdot10^{-3}$ & $4\cdot10^{-3}$ & $2\cdot
10^{-4}$ & $6\cdot10^{-4}$\\
& $10^{-6}$ & $10^{-3}$ & $4\cdot10^{-4}$ & $4\cdot10^{-4}$ & $2\cdot
10^{-5}$ & $6\cdot10^{-5}$\\
& $10^{-7}$ & $10^{-4}$ & $4\cdot10^{-5}$ & $4\cdot10^{-5}$ & $2\cdot
10^{-6}$ & $6\cdot10^{-6}$ \\[3pt]
$-3$ 
& $10^{-5}$ & $3\cdot10^{-4}$ & $1\cdot10^{-4}$ & $9\cdot10^{-5}$ &
$4\cdot10^{-6}$ & $2\cdot10^{-5}$\\
& $10^{-6}$ & $8\cdot10^{-6}$ & $3\cdot10^{-6}$ & $3\cdot10^{-6}$ &
$5\cdot10^{-8}$ & $5\cdot10^{-7}$ \\
& $10^{-7}$ & $7\cdot10^{-7}$ & $3\cdot10^{-7}$ & $2\cdot10^{-7}$ &
$1\cdot10^{-8}$ & $4\cdot10^{-8}$\\
\hline
\end{tabular*}
\end{table}

\subsection{Models with four or six parameters}\label{subsec:app_sub3}
Demidenko (\citeyear{d04}) used a double exponential model to characterize the
regrowth of tumor after radiation. The natural logarithm of tumor
volume can be modeled using a nonlinear regression model with mean
\[
\eta(x,\bftheta)=\theta_1+\log\bigl(\theta_2e^{\theta_3x}+(1-
\theta _2)e^{-\theta_4x}\bigr),
\]
where $0\le x\in[L,U]$ is the time, $\theta_1$ is the logarithm of the
initial tumor volume, $0<\theta_2<1$ is the proportional contribution
of the first compartment, and $\theta_3, \theta_4>0$ are cell
proliferation and death rates.

\citet{d06} used the LINEXP model to characterize tumor growth
delay and regrowth. The model was described in Section~\ref{subsec:locally_sub1} and re-presented below:
\[
\eta(x,\bftheta)=\theta_1+\theta_2e^{\theta_3x}+
\theta_4x.
\]

\citet{lb11} considered $D$- and $c$-optimal designs for
these two models, but our approach yields more general results. For
both models, Theorem~\ref{thmm:thmm1}(d) can be applied, and a complete
class consists of designs with at most four design points including
both endpoints [see \citet{ys12}]. Thus, Corollary~\ref
{corl:corl1} can again be applied on the design space $[L,U]$, and $\Phi
_p$-optimal designs for $\bftheta$ and certain $c$-optimal designs can
be found by solving for the critical points. In particular, $\bff_{-3}$
and $\bff_{-4}$ are Chebyshev systems under both models [see \citet{lb11}], thus $\bfe_3$- and $\bfe_4$-optimal designs for
both models can be found by solving for the critical points.

There is no explicit solution for the optimal designs, but numerical
solutions can be easily found using Newton's algorithm. Here, we give
some $A$-optimal designs for the LINEXP model in Table~\ref{Tab:tab4}
(the optimal designs for the LINEXP model do not depend on $\theta_1$
and $\theta_4$ since they are not involved in the information matrix).
For $D$- and $c$-optimality, our approach gives the same results as in
\citet{lb11}.

\begin{table}
\caption{$A$-optimal designs for the LINEXP model on $[0,1]$}\label{Tab:tab4}
\begin{tabular*}{\textwidth}{@{\extracolsep{\fill}}lccc@{}}
\hline
${\theta_2}$ & ${\theta_3}$ & ${(x_1, x_2, x_3, x_4)}$ &
\multicolumn{1}{c@{}}{${(\omega_1, \omega
_2, \omega_3, \omega_4)}$} \\
\hline
0.5 & $-$1 & (0, 0.220, 0.717, 1) & (0.156, 0.324, 0.344, 0.176) \\
1 & $-$1 & (0, 0.220, 0.717, 1) & (0.151, 0.319, 0.349, 0.181) \\
1 & $-$2 & (0, 0.195, 0.681, 1) & (0.146, 0.315, 0.355, 0.184) \\
\hline
\end{tabular*}
\end{table}

Consider one more example. \citet{dmw06} studied
$D$-optimal designs for exponential regression models, which are
nonlinear regression models with mean
%
\begin{equation}
\label{mod:expreg} \eta(x,\bftheta)=\sum_{s=1}^{S}
\theta_{2s-1}e^{-\theta_{2s}x},\qquad 0\le x\in[L, U],
\end{equation}
where $\theta_{2s-1}\neq0, s=1,\ldots,S, 0<\theta_2<\cdots<\theta
_{2S}$. When $S=2$ and $\theta_4/\theta_2<61.98$ or $S=3, 2\theta
_4=\theta_2+\theta_6$ and $\theta_4/\theta_2<23.72$, Theorem~\ref
{thmm:thmm1}(b) can be applied, and a complete class consists of
designs with at most $2S$ design points including the lower endpoint
$L$ [see \citet{ys12}, Theorems~3 and~4]. Moreover,
it is easy to see that the information matrix $\bfM_x$ goes to 0 when
$x$ approaches infinity, thus Corollary~\ref{corl:corl1} can be applied
on any design region $[L,\infty)$. Table~\ref{Tab:tab5} gives some
$A$-optimal designs for $S=2$.

\begin{table}
\tabcolsep=0pt
\caption{$A$- and $\bfe_2$-optimal designs for exponential regression
model on $[0,\infty)$ when $S=2,\theta_1=\theta_2=1$}\label{Tab:tab5}
\begin{tabular*}{\textwidth}{@{\extracolsep{\fill}}lcccc@{}}
\hline
{Criterion} & ${\theta_3}$ & ${\theta_4}$ &
${(x_1, x_2, x_3, x_4)}$ & ${(\omega
_1, \omega_2, \omega_3, \omega_4)}$ \\
\hline
$A$-optimality & 1 & 2 & $(0, 0.275, 1.196, 3.416)$ & (0.078, 0.178,
0.251, 0.493) \\
& 1 & 4 & $(0, 0.170, 0.768, 2.472)$ & (0.118, 0.261, 0.287, 0.334) \\
& 3 & 4 & $(0, 0.172, 0.760, 2.450)$ & (0.083, 0.199, 0.296, 0.422) \\[3pt]
$\bfe_2$-optimality & 1 & 2 & $(0, 0.273, 1.197, 3.425)$ & (0.054,
0.124, 0.200, 0.623)\\
& 1 & 4 & $(0, 0.168, 0.769, 2.492)$ & (0.033, 0.082, 0.201, 0.683)\\
& 3 & 4 & $(0, 0.168, 0.769, 2.492)$ & (0.033, 0.082, 0.201, 0.683)\\
\hline
\end{tabular*}
\end{table}

For $c$-optimality, first we have
\begin{eqnarray*}
&&\bff(x,\bftheta)\\
&&\qquad=\cases{ %
\bigl(e^{-\theta_2x},-
\theta_1xe^{-\theta_2x},e^{-\theta_4x},-\theta
_3xe^{-\theta_4x}\bigr), & \quad $S=2,$
\vspace*{2pt}\cr
\bigl(e^{-\theta_2x},-\theta_1xe^{-\theta_2x},e^{-\theta_4x},-
\theta _3xe^{-\theta_4x}, e^{-\theta_6x},-\theta_5xe^{-\theta_6x}
\bigr), & \quad $S=3.$}
\end{eqnarray*}
Both are Chebyshev systems. In addition, we can show that $\bff_{-2s},
s=1,\ldots, S$ are Chebyshev systems for $S=2$ and $S=3$, so the
$c$-optimal designs for $\theta_{2s}, s=1,\ldots, S$ on $[L,\infty)$
can be found by solving for the critical points. Table~\ref{Tab:tab5}
gives some $\bfe_2$-optimal designs for $S=2$.

Moreover, the $\Phi_p$-optimal designs for $\bftheta$ and $c$-optimal
design for $\theta_{2s}$'s are unique by Theorem~\ref{thmm:thmm4}. For
a finite design region, Theorem~\ref{thmm:thmm5} can be applied. For
example, the $A$-optimal design for $\bftheta=(1,1,1,2)^T$ on $[0,U]$
when $U\ge3.416$ is the same as in Table~\ref{Tab:tab5}; when
$U<3.416$, the optimal design is supported at 4 design points including
both $0$ and $U$.

\subsection{Polynomial regression model with $d$ parameters}\label
{subsec:app_sub4}
Yang (\citeyear{y10}) considered the general $(d-1)$th degree polynomial
regression model $P_{d-1}$ with variance $\sigma^2/\lambda(x)$ and mean
%
\begin{equation}
\label{equ:modpoly} \eta(x,\bftheta)=\theta_1+\sum
_{i=2}^{d} \theta_ix^{i-1}.
\end{equation}
For different choices of the efficiency function
$\lambda(x)$, Theorem~\ref{thmm:thmm1} gives the following complete class results
[see \citet{y10}, Theorem~9]:
\begin{longlist}[(a)]
\item[(a)] When (i) $\lambda(x)=1-x, x\in[-1,1]$ or (ii) $\lambda
(x)=e^{-x}, x\in[0,\infty)$, a complete class consists of designs with
at most $d$ design points including the left endpoint. Moreover, the
information matrix $\bfM_U(\bftheta)=0$.
\item[(b)] When $\lambda(x)=1+x$, $x\in[-1,1]$, a complete class
consists of designs with at most $d$ design points including the right
endpoint. Moreover, the information matrix $\bfM_L(\bftheta)=0$.
\item[(c)] When (i) $\lambda(x)=(1-x)^{u+1}(1+x)^{v+1}, x\in[-1,1],
u+1>0, v+1>0$ or (ii) $\lambda(x)=x^{u+1}e^{-x}, x\in[0,\infty),
u+1>0$ or (iii) $\lambda(x)=e^{-x^2}, x\in(-\infty,\infty)$ or (iv)
$\lambda(x)=(1+x^2)^{-t}, x\in(-\infty,\infty), d\le t$, a complete
class consists of designs with at most $d$ design points. Moreover, the
information matrices $\bfM_L(\bftheta)=\bfM_U(\bftheta)=0$.
\item[(d)] When $\lambda(x)\equiv1, x\in[L,U]$, a complete class
consists of designs with at most $d$ design points including both endpoints.
\end{longlist}
Corollary~\ref{corl:corl1} can be applied to the above models on the
respective (full) design regions, thus $\Phi_p$-optimal designs for
$\bftheta$ and $c$-optimal designs for $\theta_d$ can be found by
solving for the critical points. Furthermore, those designs are unique,
so Theorem~\ref{thmm:thmm5} can be used when the design regions are small.

Finally, we apply our theorems to more general optimality criteria.
\citet{ds95} considered optimal designs under nested
polynomial regression models. To be specific, suppose the degree of the
polynomial regression model is an unknown integer between 1 and $d-1$.
The $D$-optimal design $\xi_D^{\ell}$ under a given model $P_{\ell}$,
$1\le\ell\le d-1$, may not be efficient under another model with a
different degree. To take this uncertainty into consideration, the
authors proposed the following weighted optimality criteria $\Phi
_{p',\bfbeta}$:
%
\begin{equation}
\label{def:mixmodcrit} \Phi_{p',\bfbeta}(\bfM_{\xi})= \Biggl[\sum
_{\ell=1}^{d-1} \beta_{\ell
}\bigl(
\operatorname{eff}_D^{\ell}(\xi)\bigr)^{p'}
\Biggr]^{1/p'},
\end{equation}
where $p'\in[-\infty, 1]$, $\bfbeta=\{\beta_1,\ldots,\beta_{d-1}\}$ is
a prior on the set $\{1,\ldots,d-1\}$ with $\beta_{d-1}>0$,
\[
\operatorname{eff}_D^{\ell}(\xi)= \biggl(
\frac{\operatorname{det}\bfM_{\xi}^{\ell}}{\operatorname
{det}\bfM^{\ell}_{\xi_D^{\ell}}} \biggr)^{{1}/{(\ell+1)}},\qquad \ell =1,\ldots,d-1,
\]
$\bfM_{\xi}^{\ell}$ is the information matrix of $\xi$ under model
$P_{\ell}$, and $\operatorname{eff}_D^{\ell}(\xi)$ is the $D$-efficiency of $\xi
$ under model $P_{\ell}$.

\citet{ds95} gave the solution of $\Phi_{p',\bfbeta
}$-optimal design for $\lambda(x)\equiv1, x\in[-1,1]$. The solution
is rather complicated, and it requires knowledge of canonical moments.
An alternative way is to use Theorem~\ref{thmm:thmm3}, and it can be
applied to more general settings.

First, the $D$-efficiency in the definition of $\Phi_{p',\bfbeta}$ can
be generalized to any $\Phi_p$-efficiency, $p\in(-\infty,1]$ (e.g.,
$A$-efficiency when $p=-1$), and we denote the resulting optimality
criteria as $\Phi_{p,p',\bfbeta}$. Second, the efficiency function
$\lambda(x)$ can be generalized to any function in cases (a)$\sim$(d)
in this subsection, where $x$ belongs to the respective (full) design regions.

Under this general setting, $\Phi_{p,p',\bfbeta}$ always satisfies
Assumption~\ref{asmp:assumption A} about optimality criteria in
Section~\ref{sec:locally} [see \citet{p93}, page~285]. Moreover,
while this optimality criterion is defined on a mixture of different
models, these models are nested within the largest model $P_{d-1}$,
thus our complete class result for $P_{d-1}$ can be applied to $\Phi
_{p,p',\bfbeta}$. Finally, to use Theorem~\ref{thmm:thmm3}, any $\Phi
_{p,p',\bfbeta}$-optimal design must have at least $d$ support points.
This requirement is reasonable since otherwise the optimal design will
not be able to estimate the model $P_{d-1}$, which may be the true
model. To meet the requirement, it is sufficient to restrict ourselves
to $p,p'\in(-\infty, 0]$, since any singular matrix will result in
$\Phi_{p,p',\bfbeta}$ to be 0. So by Theorem~\ref{thmm:thmm3}, $\Phi
_{p,p',\bfbeta}$-optimal designs for models in cases (a)$\sim$(d) of
this subsection can be found by solving for the critical points. Some
examples are given in Table~\ref{Tab:tab6} for the case $\lambda
(x)=1-x^2, x\in[-1,1], p=-1$ [i.e., for $A$-efficiency in~(\ref
{def:mixmodcrit})], $d=4$ and $\bfbeta$ a uniform prior.

\begin{table}
\tabcolsep=0pt
\caption{$\Phi_{p,p',\bfbeta}$-optimal designs for polynomial
regression models}\label{Tab:tab6}
\begin{tabular*}{\textwidth}{@{\extracolsep{\fill}}lccc@{}}
\hline
${p'}$ & $(x_1, x_2, x_3, x_4)$ &
$(\omega_1, \omega_2, \omega
_3, \omega_4)$ &
$(\operatorname{eff}_A^1(\xi)$, $\operatorname{eff}_A^2(\xi)$,
$\operatorname{eff}_A^3(\xi)) $\\
\hline
0 & ($-$0.860, $-$0.346, 0.346, 0.860) & (0.263, 0.237, 0.237, 0.263)
& (0.692, 0.745, 0.902) \\ 
$-1$ & ($-$0.854, $-$0.343, 0.343, 0.854) & (0.268, 0.232, 0.232,
0.268) & (0.701, 0.753, 0.879) \\
$-3$ & ($-$0.846, $-$0.339, 0.339, 0.846) & (0.273, 0.227, 0.227, 0.273)
& (0.714, 0.759, 0.846) \\
\hline
\end{tabular*}
\end{table}

In addition, $\Phi_{p,p',\bfbeta}$-optimality is strictly isotonic and
strictly concave on PD($d$) since $\beta_{d-1}>0$ and the $\Phi
_p$-efficiency under model $P_{d-1}$ is strictly isotonic and strictly
concave on PD($d$) for $p\in(-\infty,0]$. Hence by Theorem~\ref
{thmm:thmm4}, the optimal designs are unique. However, for smaller
design regions, the optimality criterion $\Phi_{p,p',\bfbeta}$ changes
as the design region changes. For example, when $p=0$, the design $\xi
_D^{\ell}$ changes when the design region changes, which causes $\Phi
_{p,p',\bfbeta}$ to change. So the optimal design on the full design
region cannot be used to obtain the optimal design on a smaller region
as we did in Theorem~\ref{thmm:thmm5}.

\section{Computational advantages}\label{sec:computation}
Although it is not the main motivation, our method does provide
computational advantages over other algorithms, as Newton's algorithm
is well studied, easy to program and fast. For comparison, we choose
the optimal weight exchange algorithm (OWEA) proposed in \citet{ybt13}, which is among the most general and fastest
algorithms.

OWEA algorithm starts with an initial design on a grid of the design
space, then iterates between optimizing the weights for the current set
of support points and adding a new grid point to the current support
points, until the condition for optimality in general equivalence
theorem is satisfied. The computing time increases as the grid size
$\kappa$ becomes larger. So to reduce the computing time, the authors
proposed a modified algorithm. The modified algorithm starts with a
coarse grid and finds the optimal design on the coarse grid. Based on
that, the grid near the support points of the optimal design is refined
and a more accurate optimal design is found on the finer grid. We refer
to their original and modified algorithm as OWEA I and OWEA II,
respectively. All algorithms are coded using SAS IML and run on a Dell
Desktop (2.5~GHz and 4~Gb RAM). Comparisons are made for different grid
sizes, different models and under both $A$- and $D$-optimality criterion.

First, we consider the LINEXP model given in (\ref{equ:modlinexp}). The
parameters are set to be $\bftheta=(1, 0.5, -1, 1)^T$, and the design
space is $[0,1]$. Three different grid sizes, $\kappa=100,1000$ and
10,000, are used for OWEA I and II; and for OWEA II, the initial coarse
grid sizes are chosen to be 10, 100 and 100, respectively. The
computing times are shown in Table~\ref{tab:tab7}. Note the grid size
$\kappa$ is irrelevant for the speed of Newton's algorithm.

\begin{table}
\caption{Computation time (seconds) for $A$- and $D$-optimal designs
for the LINEXP model}\label{tab:tab7}
\begin{tabular*}{\textwidth}{@{\extracolsep{\fill}}lcccccc@{}}
\hline
& \multicolumn{3}{c}{${{A}}${-optimal}} &
\multicolumn{3}{c@{}}{${{D}}${-optimal}}\\[-6pt]
& \multicolumn{3}{c}{\hrulefill} &
\multicolumn{3}{c@{}}{\hrulefill}\\
$ $ & ${\kappa=100}$ & ${\kappa=1000}$ &
${\kappa=10\mbox{{,}}000}$ & ${\kappa=100}$ &
${\kappa=1000}$ & \multicolumn{1}{c@{}}{${\kappa=10\mbox{{,}}000}$}\\
\hline
Newton's & 0.08 & 0.08 & 0.08 & 0.08 & 0.08 & 0.08 \\
OWEA I & 0.19 & 0.22 & 0.63 & 0.24 & 0.37 & 1.28\\
OWEA II & 0.17 & 0.18 & 0.21 & 0.20 & 0.23 & 0.29\\
\hline
\end{tabular*}
\end{table}

From Table~\ref{tab:tab7}, we can see all three algorithms are very
efficient in finding optimal designs. Newton's algorithm is at least
twice as fast as the other two algorithms. The speed gain is more
prominent when comparing to OWEA I, especially when the grid size
$\kappa$ is large.

Second, we consider a polynomial regression model given in (\ref
{equ:modpoly}) with $d=6$ and $\lambda(x)=1-x^2, x\in[-1,1]$. It has
more parameters than the previous example so finding optimal designs
takes longer. The results are shown in Table~\ref{tab:tab8}, with a
similar conclusion as in the previous example.

\begin{table}
\caption{Computation time (seconds) for $A$- and $D$-optimal designs
for polynomial regression model}\label{tab:tab8}
\begin{tabular*}{\textwidth}{@{\extracolsep{\fill}}lcccccc@{}}
\hline
& \multicolumn{3}{c}{${{A}}${-optimal}} &
\multicolumn{3}{c@{}}{${{D}}${-optimal}}\\[-6pt]
& \multicolumn{3}{c}{\hrulefill} &
\multicolumn{3}{c@{}}{\hrulefill}\\
$ $ & ${\kappa=100}$ & ${\kappa=1000}$ &
${\kappa=10\mbox{{,}}000}$ & ${\kappa=100}$ &
${\kappa=1000}$ & \multicolumn{1}{c@{}}{${\kappa=10\mbox{{,}}000}$}\\
\hline
Newton's & 0.17 & 0.17 & 0.17 & 0.14 & 0.14 & 0.14 \\
OWEA I & 0.33 & 0.61 & 3.49 & 0.48 & 1.09 & 4.83\\
OWEA II & 0.34 & 0.48 & 0.89 & 0.44 & 0.72 & 1.35\\
\hline
\end{tabular*}
\end{table}

\section{Discussion}\label{sec:Disscusion}
In this paper, we present a general theory for finding saturated
optimal designs based on the complete class results in \citet{ys12} as well as \citet{ds13}. While we focus on locally
optimal designs, Theorem~\ref{thmm:thmm2} also applies in a multistage
design setting, and we have constructed optimal two-stage designs for
the Michaelis--Menten model using this approach. However, unlike in the
locally optimal design case, we cannot guarantee the existence of a
feasible critical point in the multistage design context, so there is
no guarantee this approach always works in that case.

For $E$-optimality, as long as the smallest eigenvalue of the
information matrix $\bfM_{\xi^*_E}$ has multiplicity 1, where $\xi^*_E$
is the $E$-optimal design, we have that $\Phi_{-\infty}$ is smooth in a
neighborhood of $\bfM_{\xi^*_E}$, and $E$-optimal designs can still be
found by solving for the critical points. For nonlinear models, we find
the smallest eigenvalue of $\bfM_{\xi^*_E}$ often does have
multiplicity 1, but we usually do not know this ahead of time. On the
other hand, we can approach $E$-optimal designs using $\Phi_p$-optimal
designs as $|p|\rightarrow\infty$. We can show whenever $|p|\ge-\log
d/\log0.95$, the $\Phi_p$-optimal design has at least 95\%
$E$-efficiency. This is not a tight bound; in practice, we find a much
smaller $|p|$ is enough.

We now point out models that cannot be accommodated. First, this occurs
when the complete class given by Theorem~\ref{thmm:thmm1} is not small
enough. For example, in \citet{dkbb10}, $D$-optimal designs for a
nonlinear model with mean $\eta(x,\bftheta)=\theta_1+\theta_2 \exp
(x/\theta_3), x\in[L,U]$ are found to be 3-point designs with both
endpoints, whereas a complete class consists of designs with at most 3
design points including only the upper endpoint as a fixed design point
[\citet{y10}, Theorem~3]. So the $D$-optimal designs are actually on
the boundary of the $\bfZ$-space, hence no feasible critical points can
be found, and the approach fails.

Second, the method fails when the model contains multiple covariates.
In general, theoretical results are very hard to obtain for
multi-covariate models, and only a couple of papers have provided some
theoretical guidance. Specific to our approach, complete class results
similar to Theorem~\ref{thmm:thmm1} are not available. The reason is
that complete class results are built upon Chebyshev systems. However,
there is no satisfactory multidimensional generalization of the
Chebyshev system yet. While \citet{yzh11} gave complete
class results for logistic and probit models with multiple covariates,
the complete classes are not derived using multidimensional Chebyshev
systems, and they are not small enough for our method to be applied.

\begin{appendix}
\section*{Appendix: Proofs}\label{sec:Proofs}
We will prove Theorems~\ref{thmm:thmm2} and~\ref{thmm:thmm4}.
Before proving Theorem~\ref{thmm:thmm2}, we first provide a lemma. This
lemma is easier stated in terms of $c$, but it can be translated into
$x$. Recall that Theorem~\ref{thmm:thmm1} gives the form of a complete
class. For any design $\xi$, we can find a design $\tilde{\xi}=\{(\tilde
{c}_j,\tilde{\omega}_j)\}_{j=1}^{m}$ in the complete class that is
noninferior ($\bfM_{\tilde{\xi}}\ge\bfM_{\xi}$).

In particular, for $\xi$ specified in Lemma~\ref{lem:lem1}, let $\Psi
_0(c)\equiv1$, a design $\tilde{\xi}$ can be found by solving the
following nonlinear equation system [see \citet{ys12} and
\citet{ds13}]:
%
\begin{equation}
\label{equ:equsystem} \sum_i \omega_i
\Psi_{\ell}(c_i)=\sum_j
\tilde{\omega}_j\Psi_{\ell
}(\tilde{c}_j),\qquad
\ell=0,1,\ldots,k-1,
\end{equation}
where $\tilde{c}_1$ and $\tilde{c}_{m}$ may be fixed to be boundary
points (see Lemma~\ref{lem:lem1}). Multiply both sides of (\ref
{equ:equsystem}) by a positive constant, the equation system still
holds, so we can remove the constraint of $\sum_i w_i=1$ for $\xi$ and
allow $\sum_i w_i$ to be any positive number in the following Lemma~\ref
{lem:lem1}; similarly for $\tilde{\xi}$ (but we still refer to them as
designs for convenience). Let $\bfX=(\bfc^T,\bfomega^T)^T$ be the
vector of all $c_i$'s and $\omega_i$'s in~$\xi$. Let $\bfS_1$ and $\bfS
_2$ be the sets of all possible vectors $\bfX$ corresponding to designs
in cases (1a)$\sim$(1d) and (2) of Lemma~\ref{lem:lem1}, respectively.
Further, let $\bfY$ be the vector of all $\tilde{c}_j$'s except those
fixed as boundary points (if any) and all $\tilde{\omega}_j$'s in
design $\tilde{\xi}$ given in the following Lemma~\ref{lem:lem1}. We
will define function $H$, $H(\bfX)=Y$, where $\bfX\in\bfS=\bfS
_1\cup\bfS_2$, and show this function is smooth on $\bfS$ under
certain conditions.

\begin{lemma}\label{lem:lem1}
Suppose one of the conditions in Theorem~\ref{thmm:thmm1} holds.
\begin{longlist}[(1a)]
\item[(1a)] If $k=2m-1$ and $\bfF(c)<0$, then for any design $\xi=\{
(c_i,\omega_i)\}_{i=1}^{m}, A<c_1<\cdots<c_{m}\le B, \omega_i>0$ for
$i\ge1$, there exists a noninferior design $\tilde{\xi}=\{(\tilde
{c}_j,\tilde{\omega}_j)\}_{j=1}^{m}$, where $\tilde{c}_{1}=A, \tilde
{\omega}_j>0$ for $j\ge1$, that solves (\ref{equ:equsystem}).
\item[(1b)] If $k=2m-1$ and $\bfF(c)>0$, then for any design $\xi=\{
(c_i,\omega_i)\}_{i=1}^{m}, A\le c_1<\cdots<c_{m}<B, \omega_i>0$ for
$i\ge1$, there exists a noninferior design $\tilde{\xi}=\{(\tilde
{c}_j,\tilde{\omega}_j)\}_{j=1}^{m}$, where $\tilde{c}_{m}=B, \tilde
{\omega}_j>0$ for $j\ge1$, that solves (\ref{equ:equsystem}).
\item[(1c)] If $k=2m$ and $\bfF(c)<0$, then for any design $\xi=\{
(c_i,\omega_i)\}_{i=1}^{m+1}, A\le c_1<\cdots<c_{m+1}\le B, \omega_i>0$
for $i\ge1$, there exists a noninferior design $\tilde{\xi}=\{(\tilde
{c}_j,\tilde{\omega}_j)\}_{j=1}^{m}$, where $\tilde{\omega}_j>0$ for
$j\ge1$, that solves (\ref{equ:equsystem}).
\item[(1d)] If $k=2m-2$ and $\bfF(c)>0$, then for any design $\xi=\{
(c_i,\omega_i)\}_{i=1}^{m-1}, A<c_1<\cdots<c_{m-1}<B, \omega_i>0$ for
$i\ge1$, there exists a noninferior design $\tilde{\xi}=\{(\tilde
{c}_j,\tilde{\omega}_j)\}_{j=1}^{m}$, where $\tilde{c}_{1}=A, \tilde
{c}_{m}=B,\tilde{\omega}_j>0$ for $j\ge1$, that solves (\ref{equ:equsystem}).
\end{longlist}
Such solution is unique under each case, hence $H$ is well defined on
$\bfS_1$.\vadjust{\goodbreak}
\begin{longlist}[(2)]
\item[(2)] For each case of \textup{(1a)}$\sim$\textup{(1d)}, let $\xi$ be similarly
defined as above except that there is exactly one 0 weight and all
other weights are positive. Then rewriting $\xi$ in the form of $\tilde
{\xi}$ in each corresponding case solves (\ref{equ:equsystem}) and
defines $H$ on $\bfS_2$. Moreover, $H$ is smooth on $\bfS=\bfS_1\cup
\bfS_2$.
\end{longlist}
\end{lemma}

\begin{pf}
We only prove for case (a), others being similar. First, let us
consider~(1a). From Lemma~1 in \citet{y10} [see also \citet{ds13},
Theorem~3.1], we know that a solution to (\ref
{equ:equsystem}) exists with $\tilde{c}_{1}=A, \tilde{\omega}_j>0, j\ge
1$. Moreover, $\bfF(c)<0$ implies that $\{\Psi_0, \Psi_1,\ldots,\Psi
_{2m-2}\}$ is a Chebyshev system [see \citet{ys12},
Proposition~4], thus such solution is unique. So $H$ is well defined on
$\bfS_1$. Now we show the smoothness on $\bfS_1$.

We have $\bfX=(c_1,\ldots,c_{m}, \omega_1,\ldots,\omega_{m})^T, \bfY
=(\tilde{c}_2,\ldots,\tilde{c}_{m}, \tilde{\omega}_1,\ldots,\tilde
{\omega}_{m})^T$ by definition ($\tilde{c}_1$ is excluded in $\bfY$
since it is fixed to be $A$). Subtract the left-hand side from the
right-hand side in (\ref{equ:equsystem}), we get an equation system
$G(\bfX, \bfY)=0$, where $G$ is smooth. So $\bfY=H(\bfX)$ is the
implicit function defined by $G(\bfX, \bfY)=0$. By implicit function
theorem, to ensure $H$ to be smooth, we only need the Jacobian matrix
$G_{\bfY}(\bfX,\bfY)=\partial G(\bfX, \bfY)/\partial\bfY$ to be
nonsingular, that is,
\begin{eqnarray*}
\label{equ:detFy} &&\operatorname{det}G_{\bfY}(\bfX,\bfY)
\\
&&\qquad= %
\left|\matrix{ 0 & \cdots& 0 &1 &\cdots & 1\vspace*{2pt}
\cr
\tilde{
\omega}_2\Psi_1'(\tilde{c}_2) &
\cdots& \tilde{\omega}_{m}\Psi _1'(
\tilde{c}_{m}) & \Psi_1(A) &\cdots & \Psi_1(
\tilde{c}_{m})\vspace*{2pt}
\cr
\vdots & \ddots& \vdots& \vdots & \vdots
&\ddots\vspace*{2pt}
\cr
\tilde{\omega}_2\Psi_{2m-2}'(
\tilde{c}_2) & \cdots& \tilde{\omega }_{m}
\Psi_{2m-2}'(\tilde{c}_{m}) &
\Psi_{2m-2}(A)& \cdots& \Psi _{2m-2}(\tilde{c}_{m})}
\right| %
\\
&&\qquad=\Biggl(\prod_{j=2}^{m}
\tilde{w}_j\Biggr) d(\tilde{\bfc})\neq0,
\end{eqnarray*}
where
%
\begin{equation}
\label{equ:mat} d(\tilde{\bfc})= %
\left|\matrix{ 1 &\cdots & 1 & 0 & \cdots&
0 \vspace*{2pt}
\cr
\Psi_1(A) &\cdots & \Psi_1(
\tilde{c}_{m}) & \Psi_1'(
\tilde{c}_2) & \cdots& \Psi_1'(
\tilde{c}_{m}) \vspace*{2pt}
\cr
\vdots &\ddots & \vdots& \vdots &
\ddots& \vdots\vspace*{2pt}
\cr
\Psi_{2m-2}(A)& \cdots&
\Psi_{2m-2}(\tilde{c}_{m}) & \Psi _{2m-2}'(
\tilde{c}_2) & \cdots& \Psi_{2m-2}'(
\tilde{c}_{m})} \right| %
.\hspace*{-20pt}
\end{equation}

Since $\tilde{w}_j>0$ for all $1\le j\le m$, we only need to show
$d(\tilde{\bfc})\neq0$. We first do some column manipulations to the
matrix in (\ref{equ:mat}). Subtract the first column from the second to
the $m$th column, then for the resulting matrix, subtract the second
column from the third to the $m$th column, continue doing this until
finally subtract the $(m-1)$th column from the $m$th column. Because
the determinant does not change during this process,
%
\begin{equation}\label{equ:detUc}\qquad
{\fontsize{10.5}{12.5}\selectfont{ d(\tilde{\bfc})= %
\left|\matrix{ \Psi_1(
\tilde{c}_2)-\Psi_1(A) &\cdots& \Psi_1(
\tilde{c}_{m})-\Psi _1(\tilde{c}_{m-1}) &
\vspace*{2pt}
\cr
\vdots & \ddots& \vdots& \bfD\vspace*{2pt}
\cr
\Psi_{2m-2}(\tilde{c}_2)-\Psi_{2m-2}(A) & \cdots&
\Psi_{2m-2}(\tilde {c}_{m})-\Psi_{2m-2}(
\tilde{c}_{m-1}) & } \right| %
,}}
\end{equation}
where $\bfD$ is the $(2m-2)\times(m-1)$ matrix,
\[
\bfD= %
\pmatrix{ \Psi_1'(
\tilde{c}_2) & \cdots& \Psi_1'(
\tilde{c}_{m}) \vspace*{2pt}
\cr
\vdots & \ddots& \vdots\vspace*{2pt}
\cr
\Psi_{2m-2}'(\tilde{c}_2) & \cdots&
\Psi_{2m-2}'(\tilde{c}_{m})} %
.
\]

Treat $A$ in the first column of the matrix in (\ref{equ:detUc}) as a
variable and fix everything else, then the determinant becomes a
real-valued function of $A$. Using the mean value theorem, we get
%
\begin{eqnarray}
\label{equ:detUc2}
&&
d(\tilde{\bfc})=(\tilde{c}_2-A)
\nonumber\\
&&\qquad{} \times
\left|\matrix{ \Psi_1'(\hat{c}_1)
& \Psi_1(\tilde{c}_{3})-\Psi_1(
\tilde{c}_{2}) & \cdots  \vspace*{2pt}
\cr
\vdots& \vdots&
\ddots\vspace*{2pt}
\cr
\Psi_{2m-2}'(
\hat{c}_1) & \Psi_{2m-2}(\tilde{c}_{3})-
\Psi_{2m-2}(\tilde {c}_{2}) & \cdots} \right.\\
&&\hspace*{85pt}\left.\matrix{\Psi_1(\tilde{c}_{m})-
\Psi_1(\tilde{c}_{m-1}) &\vspace*{2pt}\cr
\vdots& \bfD\vspace*{2pt}\cr
\Psi_{2m-2}(
\tilde{c}_{m})-\Psi_{2m-2}(\tilde {c}_{m-1})&}\right|
,\nonumber
\end{eqnarray}
where $A<\hat{c}_1<\tilde{c}_2$. Let $\varepsilon=\operatorname{sign }d(\tilde
{\bfc})$ be the sign of $d(\tilde{\bfc})$, treat $\tilde{c}_2$ in the
second column of the matrix in (\ref{equ:detUc2}) as a variable, and
use the mean value theorem again to obtain
\[
\varepsilon=\operatorname{sign } 
\left|\matrix{
\Psi_1'(\hat{c}_1) & \Psi_1'(
\hat{c}_2) & \cdots& \Psi_1(\tilde {c}_{m})-
\Psi_1(\tilde{c}_{m-1}) &\vspace*{2pt}
\cr
\vdots& \vdots&
\ddots& \vdots& \bfD\vspace*{2pt}
\cr
\Psi_{2m-2}'(
\hat{c}_1) & \Psi_{2m-2}'(
\hat{c}_2) & \cdots& \Psi _{2m-2}(\tilde{c}_{m})-
\Psi_{2m-2}(\tilde{c}_{m-1}) & }\right| %
,
\]
where $\tilde{c}_2<\hat{c}_2<\tilde{c}_3$. Keep on doing this, and
finally get
\[
\varepsilon=\operatorname{sign } %
\left|\matrix{ \Psi_1'(
\hat{c}_1) & \cdots& \Psi_1'(
\hat{c}_{m-1}) & \Psi_1'(\tilde
{c}_2) & \cdots& \Psi_1'(
\tilde{c}_{m})\vspace*{2pt}
\cr
\vdots & \ddots& \cdots & \cdots &
\ddots& \cdots\vspace*{2pt}
\cr
\Psi_{2m-2}'(
\hat{c}_1) & \cdots& \Psi_{2m-2}'(
\hat{c}_{m-1}) & \Psi _{2m-2}'(
\tilde{c}_2) & \cdots& \Psi_{2m-2}'(
\tilde{c}_{m})} \right| %
,
\]
and $A=\tilde{c}_1<\hat{c}_1<\tilde{c}_2<\hat{c}_2<\cdots<\hat
{c}_{m-1}<\tilde{c}_{m}$. Since $\{\Psi_1',\ldots,\Psi_{2m-2}'\}$ is a
Chebyshev system, $\varepsilon\neq0$. Hence, the Jacobian matrix is
invertible, and the function $H$ is smooth on $\bfS_1$.

Turning to case (2), without loss of generality, assume $\omega_1=0,
\omega_i>0$ for $i\ge2$. If we can show the function $H(\bfX)$ is
continuous on $\bfS_2$ and its partial derivatives can be extended
continuously to $\bfS_2$, then it can be proved that $H(\bfX)$ is also
differentiable on $\bfS_2$. So first, we prove its continuity.

To show this, for any sequence $\bfX^n=(c_1^n,\ldots,c_{m}^n,\omega
_1^n,\ldots,\omega_{m}^n)^T, n\ge1$, $\bfomega^n>0$ and $\bfX^n$
approaching $\bfX^0=(c_1,\ldots,c_{m},0,\omega_2, \ldots,\omega_{m})^T$,
we need to show $\bfY^n=(\tilde{c}_2^n,\ldots,\tilde{c}_{m}^n,\tilde{\omega
}_1^n,\ldots,\tilde{\omega}_{m}^n)^T$ approaches $\bfY
^0=(c_2,\ldots, c_{m},0,\omega_2,\ldots,\omega_{m})^T$.

By definition, we have
%
\begin{equation}
\label{equ:continuity} \sum_{i=1}^{m}
\omega_i^n\Psi_{\ell}\bigl(c_i^n
\bigr)=\sum_{j=1}^{m} \tilde{\omega
}_j^n\Psi_{\ell}\bigl(\tilde{c}_j^n
\bigr),\qquad \ell=0,\ldots,2m-2.
\end{equation}
Suppose we have $\bfY_{j_1}^n$ does not converge to $\bfY_{j_1}^0$ for
some $j_1$, then because $\bfY^n$ is a bounded sequence, there exists a
subsequence $\{n_t| t=1,2,\ldots\}$ such that $\bfY^{n_t}$ converges to
some $\bar{\bfY}^0=(\bar{c}_2,\ldots,\bar{c}_{m}, \bar{\omega}_1,\ldots,\bar
{\omega}_{m})^T$ and $\bar{\bfY}_{j_1}^0 \neq\bfY_{j_1}^0$.

Now let $n_t\rightarrow\infty$, take the limit of (\ref
{equ:continuity}) on both sides, we get
%
\begin{equation}
\label{equ:continuity2} \sum_{i=2}^{m}
\omega_i\Psi_{\ell}(c_i)=\sum
_{j=1}^{m} \bar{\omega}_j\Psi
_{\ell}(\bar{c}_j),\qquad \ell=0,\ldots,2m-2.
\end{equation}
Since $\{\Psi_0,\ldots,\Psi_{2m-2}\}$ is a Chebyshev system and the
maximum number of different support points in (\ref{equ:continuity2})
is $2m-1$, (\ref{equ:continuity2}) only holds if $\bar{\omega}_1=0, \bar
{\omega}_i=\omega_{i}, \bar{c}_i=c_{i}$ for $i\ge2$, which means $\bar
{\bfY}^0=\bfY^0$, leading to a contradiction.

Next, we show the partial derivatives can be extended continuously to
$\bfS_2$. Using the implicit function theorem, we know
\begin{eqnarray*}
\frac{\partial H(\bfX)}{\partial\bfX}&=&-G_{\bfY}^{-1}\bigl(\bfX,H(\bfX )
\bigr)G_{\bfX}\bigl(\bfX,H(\bfX)\bigr),\\
 G_{\bfX}(\bfX,\bfY)&=&
\frac{\partial
G(\bfX,\bfY)}{\partial\bfX},
\end{eqnarray*}
for $\bfX\in\bfS_1$. When $\bfX\rightarrow\bfX^0$, $H(\bfX
)\rightarrow H(\bfX^0)$ by continuity, hence $G_{\bfY}(\bfX,\break  H(\bfX
))\rightarrow G_{\bfY}(\bfX^0,H(\bfX^0))$ since $G_{\bfY}(\bfX,\bfY)$
is continuous. Furthermore,\break $G_{\bfY}(\bfX^0,H(\bfX^0))$ is nonsingular
by the similar argument as previously, therefore, $G_{\bfY}^{-1}(\bfX
,H(\bfX))\rightarrow G_{\bfY}^{-1}(\bfX^0,H(\bfX^0))$. It is easy to
see $G_{\bfX}(\bfX,H(\bfX))\rightarrow G_{\bfX}(\bfX^0,H(\bfX^0))$,
therefore, the derivative $\partial H(\bfX)/\partial\bfX\rightarrow
-G_{\bfY}^{-1}(\bfX^0,\break H(\bfX^0))\times G_{\bfX}(\bfX^0,H(\bfX^0))$, that is,
the derivative can be extended continuously to $\bfS_2$. So $H(\bfX)$
is differentiable on $\bfS_2$ and the partial derivatives are continuous.
\end{pf}

Now we are ready to prove Theorem~\ref{thmm:thmm2}; the proof is stated
in terms of $x$ to be consistent with the theorem.

\begin{pf*}{Proof of Theorem~\ref{thmm:thmm2}}
We only prove the case where the complete class consists of designs
with at most $m$ points including $L$, other cases being similar.
Assume the design $\xi^c$ given by a feasible critical point is not an
optimal design, and an optimal design exists as $\xi^{*}=\{(L, 1-\sum_{i=2}^m \omega_i^*), \{(x_i^*,\omega_i^*)\}_{i=2}^{m}\}$, where
$L<x_2^*<\cdots<x_m^*$ is a strictly increasing sequence (some of the
weights $\omega_i^*$ may be 0 if the support size of $\xi^*$ is less
than $m$). We have $\Phi(\bfM_{\xi^{*}})>\Phi(\bfM_{\xi^{c}})$.
Consider the linear combination of the two designs, $\xi_{\epsilon
}=\epsilon\xi^{*}+(1-\epsilon)\xi^{c}$, $0 \le\epsilon\le1$, so
\[
\xi_{\epsilon}=\Biggl\{\Biggl(L, 1- (1-\epsilon)\sum
_{i=2}^m \omega_i^{c}-
\epsilon \sum_{i=2}^m
\omega_i^{*}\Biggr),\bigl\{\bigl(x_i^{c},(1-
\epsilon)\omega_i^{c}\bigr)\bigr\} _{i=2}^m,
\bigl\{\bigl(x_i^{*},\epsilon\omega_i^{*}
\bigr)\bigr\}_{i=2}^m\Biggr\}.
\]
By the concavity of the optimality criterion $\Phi$, we have
%
\begin{equation}
\label{concave} \Phi(\bfM_{\xi_{\epsilon}}) \ge(1-\epsilon)\Phi(
\bfM_{\xi
^{c}})+\epsilon\Phi(\bfM_{\xi^{*}}).
\end{equation}
Utilizing (\ref{concave}), we can get
%
\begin{equation}
\label{ge0} \frac{\Phi(\bfM_{\xi_{\epsilon}})-\Phi(\bfM_{\xi^{c}})}{\epsilon} \ge \Phi(\bfM_{\xi^{*}})-\Phi(
\bfM_{\xi^{c}}) >0.
\end{equation}
Now, if we can find a series of designs with $m$ support points, $\tilde
{\xi}_{\epsilon}=\{(L, 1-\sum_{i=2}^m \omega_{i,\epsilon}), \{
(x_{i,\epsilon}, \omega_{i,\epsilon})\}_{i=2}^m\}$, $\epsilon\ge0$
belongs to a neighborhood of 0, such that:
\begin{longlist}[1.]
\item[1.]$\Phi(\bfM_{\tilde{\xi}_{\epsilon}})\ge\Phi(\bfM_{\xi_{\epsilon}})$;
\item[2.]$\bfZ_{\epsilon}=(\bfx_{\epsilon},\bfomega_{\epsilon})$ depends
smoothly on $\epsilon$, where $\bfx_{\epsilon}=(x_{2,\epsilon},\ldots
,x_{m,\epsilon})$, $\bfomega_{\epsilon}=(\omega_{2,\epsilon},\ldots
,\omega_{m,\epsilon})$;
\item[3.]$\bfZ_0=\bfZ^c=(\bfx^c,\bfomega^c)$, thus $\tilde{\xi}_0=\xi^c$.
\end{longlist}
Then, applying (\ref{ge0}), we obtain
\[
\frac{\Phi(\bfM_{\tilde{\xi}_{\epsilon}})-\Phi(\bfM_{\tilde{\xi
}_{0}})}{\epsilon} \ge\frac{\Phi(\bfM_{\xi_{\epsilon}})-\Phi(\bfM_{\xi
^{c}})}{\epsilon} \ge \Phi(\bfM_{\xi^{*}})-\Phi(
\bfM_{\xi^{c}}) >0.
\]
Because $\tilde{\xi}_{\epsilon}$ has $m\ge d$ support points, $\bfM
_{\tilde{\xi}_{\epsilon}}$ must belong to PD($d$). By our smoothness
assumption of $\Phi$, $\Phi(\bfM_{\tilde{\xi}_{\epsilon}})$ is a smooth
function of $\epsilon$. Take the limit as $\epsilon\rightarrow0$, it gives
%
\begin{equation}
\label{dge0} \frac{\partial\Phi(\bfM_{\tilde{\xi}_{\epsilon}})}{\partial\epsilon
} \bigg|_{\epsilon=0}>0.
\end{equation}
%
On the other hand, by our definition, $\Phi(\bfM_{\tilde{\xi}_{\epsilon
}})=\tilde{\Phi}(\bfZ_{\epsilon})$. Applying the chain rule and using
the fact that $\bfZ_0=\bfZ^c$ is a critical point of $\tilde{\Phi}(\bfZ
)$, we can get
\[
\frac{\partial\Phi(\bfM_{\tilde{\xi}_{\epsilon}})}{\partial\epsilon
} \bigg|_{\epsilon=0}=\frac{\partial\tilde{\Phi}(\bfZ_{\epsilon
})}{\partial\epsilon} \bigg|_{\epsilon=0}=
\frac{\partial\tilde{\Phi
}(\bfZ)}{\partial\bfZ} \bigg|_{\bfZ=\bfZ_0}\frac{\partial\bfZ_{\epsilon
}}{\partial\epsilon} \bigg|_{\epsilon=0}=0.
\]
This contradicts with (\ref{dge0}). Hence, $\xi^{c}$ must be an optimal design.

To find such designs $\tilde{\xi}_{\epsilon}$, first, if the design $\xi
^*$ does not have new design points other than those in $\xi^c$, that
is, $\forall2 \le i\le m$, we have either $\omega_i^*=0$ or $x_i^* \in
\bfx^c$, then the design $\xi_{\epsilon}$ is itself a design with $m$
support points, we can simply let $\tilde{\xi}_{\epsilon}=\xi_{\epsilon
}$, and conditions 1 $\sim$ 3 are satisfied.

Otherwise, suppose we have $r>0$ new design points
$x_{i_1}^*,\ldots,x_{i_r}^*$ introduced by~$\xi^*$, with $\omega
_{i_k}^*>0, k=1,\ldots, r$. Let $\delta_{\mathit{ii'}}=1$ if $x_i^c=x_{i'}^*$ and
$0$ otherwise. Rewrite the design $\xi_{\epsilon}$ as
\begin{eqnarray*}
\xi_{\epsilon} &=& \Biggl\{\Biggl(L, 1- (1-\epsilon)\sum
_{i=2}^m \omega _i^{c}-
\epsilon\sum_{i=2}^m
\omega_i^{*}\Biggr),\\
&&\quad{} \Biggl\{\Biggl(x_i^c,
(1-\epsilon)\omega _i^c+\epsilon\sum
_{i'=1}^m \omega_{i'}^*
\delta_{\mathit{ii'}}\Biggr)\Biggr\}_{i=2}^m\Biggr\}
\\
&&{}\cup\bigl\{\bigl(x_{i_k}^*, \epsilon\omega_{i_k}^*\bigr)
\bigr\}_{k=1}^r
\\
&=& \bigl\{\bigl(L, \omega_{1,\epsilon}^{(0)}\bigr), \bigl\{
\bigl(x_{i,\epsilon}^{(0)}, \omega _{i,\epsilon}^{(0)}
\bigr)\bigr\}_{i=2}^m \bigr\} \cup\bigl\{
\bigl(x_{i_k}^*, \epsilon\omega _{i_k}^*\bigr)\bigr
\}_{k=1}^r,
\end{eqnarray*}
where the second equation simply renames the design points and design
weights. It is easy to verify that conditions $2\sim3$ are satisfied
for $\bfZ_{\epsilon}^{(0)}=(\bfx_{\epsilon}^{(0)}, \bfomega_{\epsilon
}^{(0)})=(x_{2,\epsilon}^{(0)},\ldots,x_{m,\epsilon}^{(0)},\omega
_{2,\epsilon}^{(0)},\ldots,\omega_{m,\epsilon}^{(0)})$.

To find the desired $m$-point design $\tilde{\xi}_{\epsilon}$, we need
to reduce the number of design points in a ``smooth'' way. We reduce
one point at a time. First, consider the design $\{(x_{i,\epsilon
}^{(0)}, \omega_{i,\epsilon}^{(0)})\}_{i=2}^m \cup\{(x_{i_1}^*,
\epsilon\omega_{i_1}^*)\}$, all the weights are positive when
$0<\epsilon<1$, and when $\epsilon=0$, only one weight is 0. So
applying Lemma~\ref{lem:lem1} to this design we can get a new design $\{
(L, \omega_{1,\epsilon}^{(1)}),\{(x_{i,\epsilon}^{(1)},\omega
_{i,\epsilon}^{(1)})\}_{i=2}^m\}$ that is noninferior, and conditions
$2\sim3$ are satisfied for $\bfZ_{\epsilon}^{(1)}=(\bfx_{\epsilon
}^{(1)}, \bfomega_{\epsilon}^{(1)})$, where $\bfomega_{\epsilon
}^{(1)}>0$ for $0\le\epsilon<1$.

Next, we add point $x_{i_2}^*$ to $\{(x_{i,\epsilon}^{(1)},\omega
_{i,\epsilon}^{(1)})\}_{i=2}^m$ (we can always assume $x_{i_2}^*$ is a
new point to $\bfx_{\epsilon}^{(1)}$ by taking $\epsilon$ small
enough). Again, all the weights are positive when $\epsilon>0$, and
when $\epsilon=0$, only one weight is 0. Use the same method to reduce
one design point again. Keep on doing this until all $r$ new points
have been added and reduced, and we finally get $\tilde{\xi}_{\epsilon
}=\{(L, 1-\sum_{i=2}^m \omega_{i,\epsilon}^{(r)}), \{(x_{i,\epsilon
}^{(r)}, \omega_{i,\epsilon}^{(r)})\}_{i=2}^{m}\}$, that is not
inferior to $\xi_{\epsilon}$, with the conditions $1\sim3$ satisfied.
\end{pf*}

Finally, we prove Theorem~\ref{thmm:thmm4}, the proof is stated in
terms of $c$ for convenience.

\begin{pf*}{Proof of Theorem~\ref{thmm:thmm4}}
We only consider the case of Theorem~\ref{thmm:thmm1}(a). First, $\xi
^*$ must belong to the complete class. Otherwise, we can find a design
$\tilde{\xi}^*$ with $\bfM_{\tilde{\xi}^*}\ge\bfM_{\xi^*}$ and $\bfM
_{\tilde{\xi}^*}\neq\bfM_{\xi^*}$. Because $\xi^*$ has at least $d$
support points, $\bfM_{\xi^*}$ is positive definite. Since $\Phi$ is
strictly isotonic on PD($d$), we have $\Phi(\bfM_{\tilde{\xi}^*})>\Phi
(\bfM_{\xi^*})$, which is a contradiction.

Now suppose there is another optimal design $\tilde{\xi}^*$.

(i) If $\tilde{\xi}^*$ also has at least $d$ support points, then it
also belongs to the complete class by previous arguments, and we can
write $\xi^*=\{(c_i^*,\omega_i^*)\}_{i=1}^{m}$, $\tilde{\xi}^*=\{(\tilde
{c}_i^*,\tilde{\omega}_i^*)\}_{i=1}^{m}$, $c_1^*=\tilde{c}_1^*=A$. By
strict concavity, we must have $\bfM_{\xi^*}\propto\bfM_{\tilde{\xi
}^*}$ since otherwise $\Phi(\alpha\bfM_{\xi^*}+(1-\alpha)\bfM_{\tilde
{\xi}^*})>\alpha\Phi(\bfM_{\xi^*})+(1-\alpha)\Phi(\bfM_{\tilde{\xi
}^*})=\Phi(\bfM_{\xi^*}) \mbox{ for all } \alpha\in(0,1)$. Let $\bfM
_{\xi^*}=\delta\bfM_{\tilde{\xi}^*}$, then $\Phi(\delta\bfM_{\tilde{\xi
}^*})=\Phi(\bfM_{\tilde{\xi}^*})$. The strict isotonicity of $\Phi$
implies $\delta=1$, hence $\bfM_{\xi^*}= \bfM_{\tilde{\xi}^*}$ and $\bfC
_{\xi^*}=\bfC_{\tilde{\xi}^*}$. Then we have~(\ref{equ:equsystem})
holds. Because $\bfF(c)<0$, $\{\Psi_0,\ldots,\Psi_{2m-2}\}$ is a
Chebyshev system. The maximum number of different support points in
(\ref{equ:equsystem}) is $2m-1$, so (\ref{equ:equsystem}) only holds if
the design points and weights on two sides of the equations are equal,
which means $\xi^*=\tilde{\xi}^*$.

(ii) If $\tilde{\xi}^*$ has less than $d$ support points, let $\xi
_{\alpha}=\alpha\xi^*+(1-\alpha)\tilde{\xi}^*, 0<\alpha<1$. By
concavity, $\xi_{\alpha}$ is also an optimal design, moreover, it has
at least $d$ support points. Thus following the arguments in case (i),
we have $\xi_{\alpha}=\xi^*$, which means $\xi^*=\tilde{\xi}^*$. This
contradicts with the fact that $\tilde{\xi}^*$ has less than $d$
support points.
\end{pf*}
\end{appendix}




\printaddresses

\begin{thebibliography}{22}

\bibitem[\protect\citeauthoryear{de~la Garza}{1954}]{d54}
\begin{barticle}[mr]
\bauthor{\bsnm{de~la Garza},~\bfnm{A.}\binits{A.}}
(\byear{1954}).
\btitle{Spacing of information in polynomial regression}.
\bjournal{Ann. Math. Statist.}
\bvolume{25}
\bpages{123--130}.
\bid{issn={0003-4851}, mr={0060777}}
\end{barticle}
%
\bptok{imsref}%
\endbibitem

\bibitem[\protect\citeauthoryear{Demidenko}{2004}]{d04}
\begin{bbook}[mr]
\bauthor{\bsnm{Demidenko},~\bfnm{Eugene}\binits{E.}}
(\byear{2004}).
\btitle{Mixed Models: Theory and Applications}.
\bpublisher{Wiley},
\blocation{Hoboken, NJ}.
\bid{doi={10.1002/0471728438}, mr={2077875}}
\end{bbook}
%
\bptok{imsref}%
\endbibitem

\bibitem[\protect\citeauthoryear{Demidenko}{2006}]{d06}
\begin{barticle}[mr]
\bauthor{\bsnm{Demidenko},~\bfnm{Eugene}\binits{E.}}
(\byear{2006}).
\btitle{The assessment of tumour response to treatment}.
\bjournal{J. R. Stat. Soc. Ser. C. Appl. Stat.}
\bvolume{55}
\bpages{365--377}.
\bid{doi={10.1111/j.1467-9876.2006.00541.x}, issn={0035-9254}, mr={2224231}}
\end{barticle}
%
\bptok{imsref}%
\endbibitem

\bibitem[\protect\citeauthoryear{Dette}{1997}]{d97}
\begin{barticle}[mr]
\bauthor{\bsnm{Dette},~\bfnm{Holger}\binits{H.}}
(\byear{1997}).
\btitle{Designing experiments with respect to ``standardized'' optimality criteria}.
\bjournal{J. R. Stat. Soc. Ser. B Stat. Methodol.}
\bvolume{59}
\bpages{97--110}.
\bid{doi={10.1111/1467-9868.00056}, issn={0035-9246}, mr={1436556}}
\end{barticle}
%
\bptok{imsref}%
\endbibitem


\bibitem[\protect\citeauthoryear{Dette et~al.}{2008}]{dbpp08}
\begin{barticle}[mr]
\bauthor{\bsnm{Dette},~\bfnm{Holger}\binits{H.}},
\bauthor{\bsnm{Bretz},~\bfnm{Frank}\binits{F.}},
\bauthor{\bsnm{Pepelyshev},~\bfnm{Andrey}\binits{A.}} \AND
\bauthor{\bsnm{Pinheiro},~\bfnm{Jos{\'e}}\binits{J.}}
(\byear{2008}).
\btitle{Optimal designs for dose-finding studies}.
\bjournal{J. Amer. Statist. Assoc.}
\bvolume{103}
\bpages{1225--1237}.
\bid{doi={10.1198/016214508000000427}, issn={0162-1459}, mr={2462895}}
\end{barticle}
%
\bptok{imsref}%
\endbibitem

\bibitem[\protect\citeauthoryear{Dette et~al.}{2010}]{dkbb10}
\begin{barticle}[mr]
\bauthor{\bsnm{Dette},~\bfnm{H.}\binits{H.}},
\bauthor{\bsnm{Kiss},~\bfnm{C.}\binits{C.}},
\bauthor{\bsnm{Bevanda},~\bfnm{M.}\binits{M.}} \AND
\bauthor{\bsnm{Bretz},~\bfnm{F.}\binits{F.}}
(\byear{2010}).
\btitle{Optimal designs for the emax, log-linear and exponential models}.
\bjournal{Biometrika}
\bvolume{97}
\bpages{513--518}.
\bid{doi={10.1093/biomet/asq020}, issn={0006-3444}, mr={2650755}}
\end{barticle}
%
\bptok{imsref}%
\endbibitem

\bibitem[\protect\citeauthoryear{Dette, Melas and Wong}{2006}]{dmw06}
\begin{barticle}[mr]
\bauthor{\bsnm{Dette},~\bfnm{Holger}\binits{H.}},
\bauthor{\bsnm{Melas},~\bfnm{Viatcheslav~B.}\binits{V.~B.}} \AND
\bauthor{\bsnm{Wong},~\bfnm{Weng~Kee}\binits{W.~K.}}
(\byear{2006}).
\btitle{Locally {$D$}-optimal designs for exponential regression models}.
\bjournal{Statist. Sinica}
\bvolume{16}
\bpages{789--803}.
\bid{issn={1017-0405}, mr={2281302}}
\end{barticle}
%
\bptok{imsref}%
\endbibitem


\bibitem[\protect\citeauthoryear{Dette and Melas}{2011}]{dm11}
\begin{barticle}[mr]
\bauthor{\bsnm{Dette},~\bfnm{Holger}\binits{H.}} \AND
\bauthor{\bsnm{Melas},~\bfnm{Viatcheslav~B.}\binits{V.~B.}}
(\byear{2011}).
\btitle{A note on the de la {G}arza phenomenon for locally optimal designs}.
\bjournal{Ann. Statist.}
\bvolume{39}
\bpages{1266--1281}.
\bid{doi={10.1214/11-AOS875}, issn={0090-5364}, mr={2816354}}
\end{barticle}
%
\bptok{imsref}%
\endbibitem


\bibitem[\protect\citeauthoryear{Dette and Schorning}{2013}]{ds13}
\begin{barticle}[mr]
\bauthor{\bsnm{Dette},~\bfnm{Holger}\binits{H.}} \AND
\bauthor{\bsnm{Schorning},~\bfnm{Kirsten}\binits{K.}}
(\byear{2013}).
\btitle{Complete classes of designs for nonlinear regression models and principal representations of moment spaces}.
\bjournal{Ann. Statist.}
\bvolume{41}
\bpages{1260--1267}.
\bid{doi={10.1214/13-AOS1108}, issn={0090-5364}, mr={3113810}}
\end{barticle}
%
\bptok{imsref}%
\endbibitem

\bibitem[\protect\citeauthoryear{Dette and Studden}{1995}]{ds95}
\begin{barticle}[mr]
\bauthor{\bsnm{Dette},~\bfnm{Holger}\binits{H.}} \AND
\bauthor{\bsnm{Studden},~\bfnm{William~J.}\binits{W.~J.}}
(\byear{1995}).
\btitle{Optimal designs for polynomial regression when the degree is not known}.
\bjournal{Statist. Sinica}
\bvolume{5}
\bpages{459--473}.
\bid{issn={1017-0405}, mr={1347600}}
\end{barticle}
%
\bptok{imsref}%
\endbibitem






\bibitem[\protect\citeauthoryear{Elfving}{1952}]{e52}
\begin{barticle}[mr]
\bauthor{\bsnm{Elfving},~\bfnm{G.}\binits{G.}}
(\byear{1952}).
\btitle{Optimum allocation in linear regression theory}.
\bjournal{Ann. Math. Statist.}
\bvolume{23}
\bpages{255--262}.
\bid{issn={0003-4851}, mr={0047998}}
\end{barticle}
%
\bptok{imsref}%
\endbibitem

\bibitem[\protect\citeauthoryear{Karlin and Studden}{1966}]{ks66}
\begin{bbook}[mr]
\bauthor{\bsnm{Karlin},~\bfnm{Samuel}\binits{S.}} \AND
\bauthor{\bsnm{Studden},~\bfnm{William~J.}\binits{W.~J.}}
(\byear{1966}).
\btitle{Tchebycheff Systems: {W}ith Applications in Analysis and Statistics}.
\bseries{Pure and Applied Mathematics}
\bvolume{XV}.
\bpublisher{Interscience},
\blocation{New York}.
\bid{mr={0204922}}
\end{bbook}
%
\bptok{imsref}%
\endbibitem


\bibitem[\protect\citeauthoryear{Kiefer and Wolfowitz}{1965}]{kw65}
\begin{barticle}[mr]
\bauthor{\bsnm{Kiefer},~\bfnm{J.}\binits{J.}} \AND
\bauthor{\bsnm{Wolfowitz},~\bfnm{J.}\binits{J.}}
(\byear{1965}).
\btitle{On a theorem of {H}oel and {L}evine on extrapolation designs}.
\bjournal{Ann. Math. Statist.}
\bvolume{36}
\bpages{1627--1655}.
\bid{issn={0003-4851}, mr={0185769}}
\end{barticle}
%
\bptok{imsref}%
\endbibitem

\bibitem[\protect\citeauthoryear{Li and Balakrishnan}{2011}]{lb11}
\begin{barticle}[mr]
\bauthor{\bsnm{Li},~\bfnm{Gang}\binits{G.}} \AND
\bauthor{\bsnm{Balakrishnan},~\bfnm{N.}\binits{N.}}
(\byear{2011}).
\btitle{Optimal designs for tumor regrowth models}.
\bjournal{J. Statist. Plann. Inference}
\bvolume{141}
\bpages{644--654}.
\bid{doi={10.1016/j.jspi.2010.07.009}, issn={0378-3758}, mr={2732935}}
\end{barticle}
%
\bptok{imsref}%
\endbibitem

\bibitem[\protect\citeauthoryear{Pukelsheim}{1993}]{p93}
\begin{bbook}[mr]
\bauthor{\bsnm{Pukelsheim},~\bfnm{Friedrich}\binits{F.}}
(\byear{1993}).
\btitle{Optimal Design of Experiments}.
\bpublisher{Wiley},
\blocation{New York}.
\bid{mr={1211416}}
\end{bbook}
%
\bptok{imsref}%
\endbibitem

\bibitem[\protect\citeauthoryear{Studden}{1968}]{s68}
\begin{barticle}[mr]
\bauthor{\bsnm{Studden},~\bfnm{W.~J.}\binits{W.~J.}}
(\byear{1968}).
\btitle{Optimal designs on {T}chebycheff points}.
\bjournal{Ann. Math. Statist.}
\bvolume{39}
\bpages{1435--1447}.
\bid{issn={0003-4851}, mr={0231497}}
\end{barticle}
%
\bptok{imsref}%
\endbibitem

\bibitem[\protect\citeauthoryear{Yang}{2010}]{y10}
\begin{barticle}[mr]
\bauthor{\bsnm{Yang},~\bfnm{Min}\binits{M.}}
(\byear{2010}).
\btitle{On the de la {G}arza phenomenon}.
\bjournal{Ann. Statist.}
\bvolume{38}
\bpages{2499--2524}.
\bid{doi={10.1214/09-AOS787}, issn={0090-5364}, mr={2676896}}
\end{barticle}
%
\bptok{imsref}%
\endbibitem

\bibitem[\protect\citeauthoryear{Yang, Biedermann and Tang}{2013}]{ybt13}
\begin{barticle}[mr]
\bauthor{\bsnm{Yang},~\bfnm{Min}\binits{M.}},
\bauthor{\bsnm{Biedermann},~\bfnm{Stefanie}\binits{S.}} \AND
\bauthor{\bsnm{Tang},~\bfnm{Elina}\binits{E.}}
(\byear{2013}).
\btitle{On optimal designs for nonlinear models: A~general and efficient algorithm}.
\bjournal{J. Amer. Statist. Assoc.}
\bvolume{108}
\bpages{1411--1420}.
\bid{doi={10.1080/01621459.2013.806268}, issn={0162-1459}, mr={3174717}}
\end{barticle}
%
\bptok{imsref}%
\endbibitem

\bibitem[\protect\citeauthoryear{Yang and Stufken}{2009}]{ys09}
\begin{barticle}[mr]
\bauthor{\bsnm{Yang},~\bfnm{Min}\binits{M.}} \AND
\bauthor{\bsnm{Stufken},~\bfnm{John}\binits{J.}}
(\byear{2009}).
\btitle{Support points of locally optimal designs for nonlinear models with two parameters}.
\bjournal{Ann. Statist.}
\bvolume{37}
\bpages{518--541}.
\bid{doi={10.1214/07-AOS560}, issn={0090-5364}, mr={2488361}}
\end{barticle}
%
\bptok{imsref}%
\endbibitem

\bibitem[\protect\citeauthoryear{Yang and Stufken}{2012}]{ys12}
\begin{barticle}[mr]
\bauthor{\bsnm{Yang},~\bfnm{Min}\binits{M.}} \AND
\bauthor{\bsnm{Stufken},~\bfnm{John}\binits{J.}}
(\byear{2012}).
\btitle{Identifying locally optimal designs for nonlinear models: A~simple extension with profound consequences}.
\bjournal{Ann. Statist.}
\bvolume{40}
\bpages{1665--1681}.
\bid{doi={10.1214/12-AOS992}, issn={0090-5364}, mr={3015039}}
\end{barticle}
%
\bptok{imsref}%
\endbibitem

\bibitem[\protect\citeauthoryear{Yang, Zhang and Huang}{2011}]{yzh11}
\begin{barticle}[mr]
\bauthor{\bsnm{Yang},~\bfnm{Min}\binits{M.}},
\bauthor{\bsnm{Zhang},~\bfnm{Bin}\binits{B.}} \AND
\bauthor{\bsnm{Huang},~\bfnm{Shuguang}\binits{S.}}
(\byear{2011}).
\btitle{Optimal designs for generalized linear models with multiple design variables}.
\bjournal{Statist. Sinica}
\bvolume{21}
\bpages{1415--1430}.
\bid{doi={10.5705/ss.2009.115}, issn={1017-0405}, mr={2827529}}
\end{barticle}
%
\bptok{imsref}%
\endbibitem
\end{thebibliography}
\end{document}